\documentclass[11pt,a4paper]{smfart}
\usepackage{pinlabel}
\usepackage[francais]{babel}

\input xy
\xyoption{all}

\numberwithin{equation}{section} 

\theoremstyle{plain} 
\newtheorem{theorem}[equation]{Th{\'e}or{\`e}me} 
\newtheorem{lemma}[equation]{Lemme} 
\newtheorem{proposition}[equation]{Proposition}

\theoremstyle{definition} 
\newtheorem{example}[equation]{Exemple} 
\newtheorem{remark}[equation]{Remarque} 
\newtheorem{definition}[equation]{D{\'e}finition}

\newcommand{\gras}[1]{{\mathbb #1}} 
\newcommand{\C}{\gras{C}} 
\newcommand{\Z}{\gras{Z}} 
 
\newcommand{\R}{\gras{R}} 
\newcommand{\N}{\gras{N}}

\begin{document}

\title[Le cerf-volant d'une constellation]{Le cerf-volant d'une constellation} 
\author{\sc Patrick Popescu-Pampu} 
\address{Univ. Paris 7 Denis Diderot, Inst. de 
  Maths.-UMR CNRS 7586, {\'e}quipe "G{\'e}om{\'e}trie et dynamique" \\ 
  Site Chevaleret, Case 
  7012, 75205 Paris Cedex 13, France.} 
\email{ppopescu@math.jussieu.fr}

 
 
\date{16 Juin 2009.} 
 
\begin{abstract} 
  On consid{\`e}re un point lisse $O$ d'une surface analytique complexe $S$. 
  Une \emph{constellation} bas{\'e}e en $O$ est un ensemble de points infiniment
  voisins de $O$, centres d'une suite d'{\'e}clatements de points
  au-dessus de $O$. Les constellations finies sont cod{\'e}es habituellement de
  deux mani{\`e}res : soit {\`a} l'aide d'un \emph{diagramme d'Enriques},
  soit {\`a} l'aide du \emph{graphe dual} du diviseur obtenu en {\'e}clatant les
  points de la constellation. Il s'agit de deux arbres d{\'e}cor{\'e}s,
  codant compl{\`e}tement la combinatoire de la constellation. Des
  algorithmes de passage de l'un {\`a} l'autre sont connus, mais ils ne
  permettent pas de se repr{\'e}senter g{\'e}om{\'e}triquement leur
  relation. Nous
  associons {\`a} une constellation un complexe 
  simplicial g{\'e}om{\'e}trique de dimension deux, appel{\'e} son
  \emph{cerf-volant}, muni 
  d'une structure affine, et nous prouvons qu'il
  contient canoniquement le diagramme d'Enriques et le graphe dual. De
  plus les d{\'e}corations de ces deux arbres s'expriment  
  tr{\`e}s simplement en termes de la g{\'e}om{\'e}trie affine du
  cerf-volant. Ceci permet de comprendre g{\'e}om{\'e}triquement les
  relations entre les deux graphes, ainsi que leurs relations avec
  l'arbre valuatif de Favre et Jonsson, qui peut {\^e}tre interpr{\'e}t{\'e}
  en tant que graphe dual de la
  constellation de tous les points infiment voisins de $O$. En fait, les
  cerf-volants des 
  constellations finies se recollent en un cerf-volant infini qui est
  muni d'un feuilletage de dimension 1 dont l'espace des feuilles est
  l'arbre valuatif. La transition vers les calculs de fractions
  continues est assur{\'e}e par des plongements partiels des
  cerf-volants dans un complexe simplicial canoniquement associ{\'e} {\`a}
  une base d'un r{\'e}seau, son \emph{lotus}. Cette derni{\`e}re notion
  est bri{\`e}vement expos{\'e}e en toutes dimensions. 

\end{abstract}

\maketitle 
\pagestyle{myheadings} \markboth{{\normalsize 
P. Popescu-Pampu}} 
{{\normalsize Le cerf-volant d'une constellation}}


\section{Introduction} \label{intro} 

Depuis le travail fondateur \cite{N 75} de Max Noether, divers probl{\`e}mes de
g{\'e}om{\'e}trie birationnelle des surfaces ou des courbes planes ont {\'e}t{\'e} 
{\'e}tudi{\'e}s {\`a} l'aide de la notion de \emph{point infiniment voisin}
d'un point lisse donn{\'e} sur une surface. La premi{\`e}re {\'e}tude syst{\'e}matique de
cette notion dans un tra{\^\i}t{\'e} de g{\'e}om{\'e}trie alg{\'e}brique
semble {\^e}tre celle faite par Enriques et Chisini \cite{EC 18}. 

En g{\'e}om{\'e}trie birationnelle 
apparaissent naturellement des suites
d'{\'e}clatements centr{\'e}s en des points infiniment voisins de points
lisses de surfaces. Si tous 
ces {\'e}clatements sont effectu{\'e}s au-dessus d'un unique point $O$,
leurs centres forment ce que Campillo, Gonz{\'a}lez-Sprinberg \&
Lejeune-Jalabert \cite{CGL 92} appel{\`e}rent une \emph{constellation}
de points infiniment 
voisins de $O$.  Il est important de bien comprendre la g{\'e}om{\'e}trie
des constellations, en particulier comment les points se suivent les
uns les autres lors des processus d'{\'e}clatements successifs les
faisant appara{\^\i}tre. 

Dans \cite{EC 18} (voir aussi \cite{CA 00}) fut introduit un arbre enracin{\'e},
appel{\'e} depuis  
\emph{diagramme d'Enriques}, qui retient exactement le processus
pr{\'e}c{\'e}dent: ses sommets correspondent bijectivement aux points de
la constellation, la racine correspondnat {\`a} $O$, et l'on relie deux
sommets si l'un des points 
correspondants appara{\^\i}t en {\'e}clatant l'autre. De plus, les
ar{\^e}tes sont de deux types, soit courbes, soit droites, et dans une
suite d'ar{\^e}tes droites s'{\'e}loignant de la racine de l'arbre, on
dit {\`a} chaque pas si l'on va tout droit ou si l'on part 
transversalement. Avec
ces d{\'e}corations, le diagramme d'Enriques code compl{\`e}tement la
combinatoire de la constellation. 

Ult{\'e}rieurement fut introduit un autre diagramme codant
diff{\'e}remment la combinatoire de la constellation : le \emph{graphe dual} du
diviseur obtenu en {\'e}clatant tous les points de la
constellation (voir \cite{BK 86} et \cite{W 04}). Dans ce cas, ce sont
les sommets qui sont d{\'e}cor{\'e}s 
(par les auto-intersections des composantes irr{\'e}ductibles du diviseur
exceptionnel qui leur correspondent). 

Les descriptions pr{\'e}c{\'e}dentes montrent que les sommets des deux
graphes sont naturellement en correspondance bijective. Mais de cette
mani{\`e}re les ar{\^e}tes ne se correspondent pas. D'autre part, 
comme les deux diagrammes contiennent la m{\^e}me information, on peut
en principe passer de l'un {\`a} l'autre. Des algorithmes de passage ont
{\'e}t{\'e} d{\'e}crits, utilisant tous des calculs de fractions
continues (voir \cite{CC 05}). Mais ils ne permettent pas de \emph{penser 
g{\'e}om{\'e}triquement} le lien entre les deux graphes. 

Dans cet article je d{\'e}cris une mani{\`e}re de visualiser
simultan{\'e}ment le diagramme d'Enriques et le graphe dual. Pour
cela, j'associe {\`a} chaque constellation un complexe simplicial
g{\'e}om{\'e}trique de dimension deux, dont certains points sont distribu{\'e}s en
types. Je l'appele \emph{le cerf-volant de la constellation},
compos{\'e} de pi{\`e}ces {\'e}l{\'e}mentaires 
triangulaires - \emph{les voiles {\'e}l{\'e}mentaires} - et de segments - 
\emph{les
  cordes}, recoll{\'e}es lors d'un jeu d'assemblage dict{\'e} par le
procesus d'{\'e}clatement menant au d{\'e}ploiement complet de la
constellation par {\'e}clatements successifs. Certaines cordes se retrouvent
recoll{\'e}es {\`a} l'int{\'e}rieur des voiles {\'e}l{\'e}mentaires, les autres flottant
librement. L'union des voiles
{\'e}l{\'e}mentaires forme \emph{la voilure} du cerf-volant de la
constellation. Celle-ci peut {\^e}tre munie canoniquement d'une
structure affine recollant celles des voiles {\'e}l{\'e}mentaires. Tout ceci
est expliqu{\'e} dans les sections 2--5.

Le th{\'e}or{\`e}me principal de cet article (Th{\'e}or{\`e}me \ref{isos})
montre  que \emph{le diagramme d'Enriques et le graphe dual se
plongent naturellement dans le cerf-volant}. Plus pr{\'e}cis{\'e}ment, le
diagramme d'Enriques est isomorphe au graphe form{\'e} par les cordes et
le diagramme dual {\`a} une partie du bord de la voilure. De plus, les
deux types de 
d{\'e}corations s'interpr{\`e}tent en termes tr{\`e}s simples  {\`a} l'aide de la
structure affine du cerf-volant. En 
particulier, les ar{\^e}tes droites du diagramme d'Enriques correspondent
aux cordes 
internes, et elles vont tout droit selon la convention d'Enriques et
Chisini si et seulement
si elles forment une g{\'e}od{\'e}sique pour la structure affine de la
voilure !

Dans la Section \ref{intval}  j'explique le lien de la notion de
cerf-volant avec l'arbre valuatif de Favre et Jonsson 
\cite{FJ 04}. Plus pr{\'e}cis{\'e}ment, les voilures de toutes les
constellations finies se recollent en un complexe simplicial infini,
la voilure du \emph{firmament} de $O$, c'est-{\`a}-dire de la constellation de
tous les points infiniment voisins de $O$. Cette voilure peut {\^e}tre munie
naturellement d'un feuilletage de dimension $1$, dont l'espace des
feuilles est l'arbre valuatif. 

Jusqu'{\`a} pr{\'e}sent, le passage d'un graphe {\`a} l'autre {\'e}tait
d{\'e}crit {\`a} l'aide de fractions continues. Ces calculs peuvent eux
aussi {\^e}tre interpr{\'e}t{\'e}s {\`a} l'aide du cerf-volant. Pour cela,
j'introduis dans la Section \ref{seclot} la notion de \emph{lotus}
associ{\'e} {\`a} une base d'un r{\'e}seau bidimensionnel. Il s'agit d'un
complexe simplicial bidimensionnel infini dont la structure permet de
`voir'  les d{\'e}veloppements en fractions continues et dans lesquels se
plongent les voilures des cerf-volants. Dans la Section
\ref{lotarb} j'explique bri{\`e}vement  une extension de la notion de lotus
en dimension quelconque, ainsi que la g{\'e}n{\'e}ralisation de la notion
de fraction continue sugg{\'e}r{\'e}e par cette extension.
\medskip

{\bf Remerciements.} Je remercie Charles Favre pour ses explications
concernant l'arbre valuatif, Evelia Garc{\'\i}a Barroso pour les
longues discussions que nous avons eues au sujet du vocabulaire marin
et c{\'e}leste de cet article, Bernard Teissier pour m'avoir
sugg{\'e}r{\'e} le nom de `lotus', ainsi que Monique Lejeune-Jalabert pour ses
remarques sur une version pr{\'e}liminaire de cet article.

\section{Constellations de points infiniment voisins}

Dans cette section j'explique les notions de base de \emph{points
infiniment voisins}, de \emph{points proches}, d'\emph{astres}, de
\emph{firmament} et de \emph{constellations}.  
\medskip

Dans tout ce qui suit nous travaillerons avec des surfaces 
analytiques complexes lisses. Mais nos consid{\'e}rations se transposent sans 
aucun changement au cas des surfaces alg{\'e}briques lisses sur un corps 
alg{\'e}briquement clos, pourvu que tous les points soient interpr{\'e}t{\'e}s comme des 
points ferm{\'e}s. 

Soit $(S,O)$ un germe de surface lisse. Notons par 
$\mathcal{O}_{S,O}$ son alg{\`e}bre locale et par $F_{S,O}$ le corps des
fractions de $\mathcal{O}_{S,O}$.  Soit :
$$(S_O, E_O) \stackrel{\pi_O}{\rightarrow} (S,O)$$ 
le morphisme
d'{\'e}clatement du point $O$. Les points de la courbe exceptionnelle 
$E_O:= \pi_O^{-1}(O)$ sont appel{\'e}s \emph{les
  points infiniment voisins de $O$ {\`a} hauteur} $1$, ou bien
\emph{les points directement proches} de $O$. 

\begin{definition}
Si $(\Sigma, E)\stackrel{\pi}{\rightarrow} (S,O)$ est un morphisme
compos{\'e} 
d'{\'e}clatements de points, alors un point du diviseur exceptionnel r{\'e}duit
$E:= \pi^{-1}(O)$ est appel{\'e} 
\textbf{un point infiniment voisin} de $O$. 
\end{definition}

En associant {\`a} chaque
point infiniment voisin de $O$ la valuation divisorielle de $F_{S,O}$
qui calcule la 
multiplicit{\'e} au point respectif, on d{\'e}finit  
naturellement une relation d'{\'e}quivalence sur l'ensemble des points
infiniment voisins de $O$ sur les divers {\'e}clat{\'e}s  
$(\Sigma, E)$  de $(S,O)$. Par la suite, lorsque l'on parlera de
points infiniment voisins, il s'agira d'une classe d'{\'e}quivalence de
points identifi{\'e}s de la mani{\`e}re pr{\'e}c{\'e}dente. On dira qu'une
surface  $(\Sigma, E)$ obtenue par une suite d'{\'e}clatements au-dessus
de $O$ et telle que $E$ contient un repr{\'e}sentant de la classe
d'{\'e}quivalence est \emph{un mod{\`e}le contenant le point infiniment voisin}. 

On dit qu'un point infiniment voisin de $O$ est
\emph{{\`a} hauteur $d>0$} s'il est directement proche d'un
point {\`a} hauteur $d-1$. Il est dit \emph{proche de $O$} si dans un
mod{\`e}le $\Sigma$ il
se trouve sur la transform{\'e}e stricte sur $\Sigma$ de $E_O$ (on
utilise ici le fait qu'un morphisme $\pi$ non-trivial se factorise
n{\'e}cessairement par l'{\'e}clatement $\pi_O$ de $O$). Ceci est alors vrai
pour tout mod{\`e}le le contenant. 

\begin{remark}
  La notion de point infiniment voisin a {\'e}t{\'e} introduite (sous
  l'appellation \emph{unendliche nahe einander ...}) par Max Noether dans
  \cite{N 75}, afin 
  d'{\'e}tudier les singularit{\'e}s des courbes alg{\'e}briques planes {\`a} l'aide
  de suites d'{\'e}clatements successifs. 
  Les notions de points infiniments voisins et de points proches ont
  {\'e}t{\'e} utilis{\'e}es par Enriques et Chisini \cite{EC 18} pour
  formuler des conditions de passage par des points bases pour les
  syst{\`e}mes lin{\'e}aires de courbes planes. Ult{\'e}rieurement, Zariski
  \cite{Z 38} les a reformul{\'e}es dans le langage des id{\'e}aux et les a
  utilis{\'e}es pour {\'e}tudier les id{\'e}aux primaires pour l'id{\'e}al
  maximal de $\mathcal{O}_{S,O}$. Pour une introduction {\`a} ces
  aspects on pourra consulter l'article de survol \cite{LJ 95} de
  Lejeune-Jalabert.   
\end{remark}

\begin{definition}
Notons par $\mathcal{C}_O$ l'ensemble des points infiniment voisins
de $O$, en incluant $O$ lui-m{\^e}me. Nous appelerons ses {\'e}l{\'e}ments
des \textbf{astres}, $\mathcal{C}_O$ {\'e}tant le \textbf{firmament} de
$O$. 
\end{definition}

La hauteur peut {\^e}tre vue comme une fonction :
  $$H : \mathcal{C}_{O} \rightarrow \N.$$

La relation de proximit{\'e} s'{\'e}tend naturellement {\`a} $\mathcal{C}_O$
tout entier. Chaque astre $A \in \mathcal{C}_O\setminus O$ est proche
d'\emph{un} ou de 
\emph{deux} autres astres. Afin d'{\'e}tudier les deux possibilit{\'e}s, consid{\'e}rons 
un mod{\`e}le $(\Sigma, E)\stackrel{\pi}{\rightarrow} (S,O)$ contenant 
$A$. 

$\bullet$ \emph{Si $A$ se trouve sur une seule composante irr{\'e}ductible
$E_i$ de $E$}, on l'appelle un \emph{astre libre}. Il est proche
uniquement de l'astre $A_i$ dont l'{\'e}clatement cr{\'e}e $E_i$. Nous
notons $p_D(A):= A_i$, et nous appelons ce point 
\emph{le pr{\'e}d{\'e}cesseur direct} de $A$. 

$\bullet$ \emph{Si $A$ se trouve sur deux composantes $E_i$ et $E_j$
  de $E$}, on l'appelle un \emph{astre satellite}. Dans ce cas,
$A$ est proche de deux autres astres $A_i \neq A_j$ dont les
{\'e}clatements cr{\'e}ent $E_i$ et $E_j$ respectivement. L'un d'entre
eux - supposons qu'il s'agit de $A_i$ - est n{\'e}cessairement proche
de l'autre - $A_j$. Nous notons $p_D(A):= A_i, \: p_I(A):= A_j$  et nous
appelons $A_i$ \emph{le pr{\'e}d{\'e}cesseur direct} de $A$ et $A_j$
\emph{le pr{\'e}d{\'e}cesseur indirect} de $A$. 

{\'E}tendons les d{\'e}finitions pr{\'e}c{\'e}dentes en posant
$p_D(O):=O$. Nous obtenons ainsi une application surjective :
$$p_D : \mathcal{C}_O \longrightarrow \mathcal{C}_O. $$
Elle v{\'e}rifie $H \circ p_D = H-1$ sur $\mathcal{C}_O \setminus O$ et $H(O)=0$. 
Les fibres de $p_D$ au-dessus d'un astre sont les points directement
proches de celui-ci, {\`a} l'exception de $O$, pour qui la fibre
$p_D^{-1}(O)$ contient aussi l'astre $O$ lui-m{\^e}me. 

Nous pouvons reformuler de la mani{\`e}re suivante la
d{\'e}finition d'une \emph{constellation} donn{\'e}e par Campillo,
Gonz{\'a}lez-Sprinberg et Lejeune-Jalabert dans \cite{CGL 92} et
\cite{CGL 96}: 

\begin{definition}
  Une {\bf constellation} centr{\'e}e en $O$ est un sous-ensemble
  $\mathcal{C}\subset \mathcal{C}_O$ qui est stable sous l'application $p_D$.  
\end{definition}

Par la suite nous nous restreindrons uniquement {\`a} des constellations
\emph{finies} centr{\'e}es en $O$.

\begin{remark}
  On peut penser au couple $(S,O)$ comme {\`a} une repr{\'e}sentation
  visuelle d'un point brillant sur le ciel. Jadis on interpr{\'e}tait un
  tel point comme {\'e}tant une {\'e}toile, mais avec l'av{\`e}nement des
  lunettes astronomiques puis des t{\'e}l{\'e}scopes de plus en plus
  puissants, on a appris {\`a} y voir des amas d'autres points
  brillants, pouvant {\^e}tre eux-m{\^e}mes des galaxies ou des
  {\'e}toiles. C'est pour cette raison que j'utilise le terme
  \emph{astre}, ne voulant pas pr{\'e}juger de sa nature en l'appelant
  \emph{{\'e}toile}. De plus, une constellation est une configuration
  particuli{\`e}re d'astres, ce qui montre que le vocabulaire que
  j'utilise s'adapte bien {\`a} celui introduit par Campillo,
  Gonz{\'a}lez-Sprinberg et Lejeune-Jalabert.  En fait, ces derniers
  {\'e}tudi{\`e}rent des constellations en dimension quelconque (on pourra
  consulter pour l'{\'e}tat de l'art {\`a} ce sujet le survol \cite{CGM
    09}). Je voudrais remarquer aussi qu'un vocabulaire {\`a}
  conotations c{\'e}lestes a aussi {\'e}t{\'e} utilis{\'e} par Hironaka
  \cite{H 73}. Mais sa notion d'\emph{{\'e}toile} est diff{\'e}rente de
  celle d'\emph{astre} et ce qu'il appelle \emph{vo{\^u}te {\'e}toil{\'e}e}
  est diff{\'e}rent du \emph{firmament}. 
\end{remark}

\begin{figure}[h!] 
\vspace*{6mm}
\labellist \small\hair 2pt 
  \pinlabel{$O$} at 71 438

  \pinlabel{$A_2$} at 280 382
  \pinlabel{$A_1$} at 280 450
  \pinlabel{$E_0$} at 259 333

  \pinlabel{$A_3$} at 524 460
  \pinlabel{$A_4$} at 475 431
  \pinlabel{$A_5$} at 547 431
  \pinlabel{$A_6$} at 524 370
  \pinlabel{$A_7$} at 479 370
  \pinlabel{$E_0$} at 510 340
  \pinlabel{$E_1$} at 436 445
  \pinlabel{$E_2$} at 436 383

  \pinlabel{$A_8$} at 440 210
  \pinlabel{$A_9$} at 438 163
  \pinlabel{$A_{10}$} at 509 163
  \pinlabel{$A_{11}$} at 446 96
  \pinlabel{$A_{12}$} at 417 42
  \pinlabel{$A_{13}$} at 492 42
  \pinlabel{$A_{14}$} at 536 95
  \pinlabel{$E_0$} at 502 139
  \pinlabel{$E_3$} at 557 171
  \pinlabel{$E_1$} at 452 260
  \pinlabel{$E_4$} at 378 200
  \pinlabel{$E_5$} at 378 230
  \pinlabel{$E_6$} at 566 109
  \pinlabel{$E_2$} at 467 2
  \pinlabel{$E_7$} at 388 58

  \pinlabel{$E_0$} at 159 141
  \pinlabel{$E_{10}$} at 219 172
  \pinlabel{$E_6$} at 219 109
  \pinlabel{$E_3$} at 95 200
  \pinlabel{$E_{11}$} at 138 36
  \pinlabel{$E_9$} at 34 230
  \pinlabel{$E_2$} at 53 64
  \pinlabel{$E_1$} at 71 345
  \pinlabel{$E_5$} at 124 267
  \pinlabel{$E_8$} at 124 301
  \pinlabel{$E_4$} at 12 335
  \pinlabel{$E_{12}$} at 28 43
  \pinlabel{$E_{13}$} at 28 16
  \pinlabel{$E_{14}$} at 195 36

   \pinlabel{$\pi_0$} at 170 438
  \pinlabel{$\pi_1$} at 378 438
  \pinlabel{$\pi_2$} at 536 290
  \pinlabel{$\pi_3$} at 316 140
\endlabellist 
\centering 
\includegraphics[scale=0.50]{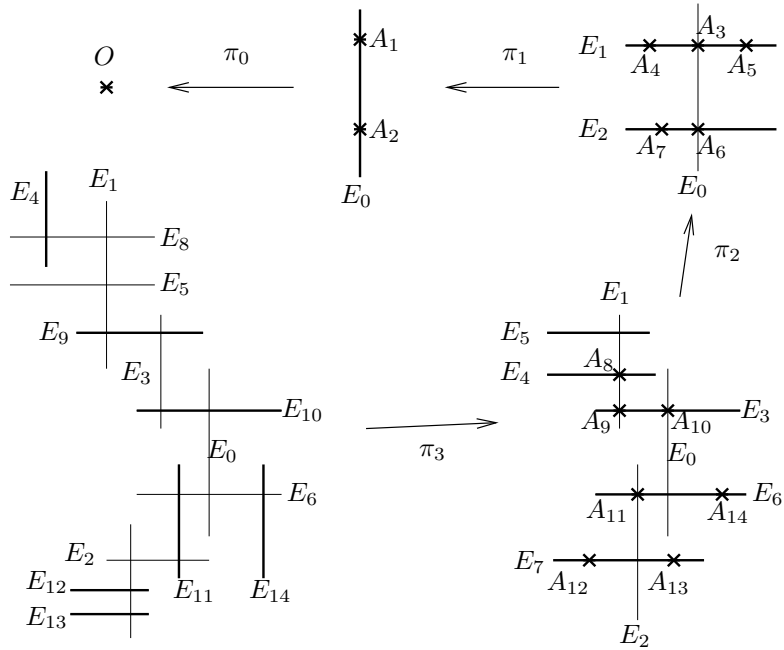} 
\vspace*{10mm} \caption{Une constellation et sa suite d'{\'e}clatements} 
\label{fig:Suitecl} 
\end{figure}

{\`A} la constellation finie $\mathcal{C}$ nous pouvons associer le 
morphisme bim{\'e}romorphe :
\begin{equation} \label{eclconst}
 \pi_{\mathcal{C}}: (S_{\mathcal{C}},
 E_{\mathcal{C}}) \rightarrow (S,O)
\end{equation}
obtenu en {\'e}clatant
successivement les astres de $\mathcal{C}$ selon leur hauteur :
on d{\'e}marre {\`a} $i=0$ et {\`a} chaque {\'e}tape $i \geq 0$ on {\'e}clate
tous les astres de $\mathcal{C}$ de hauteur $i$. Intuitivement,
cela correspond {\`a} faire des zooms successifs au 
voisinage de chaque point brillant apparaissant par le zoom ant{\'e}rieur,
pour voir si ce point correspond plut{\^o}t {\`a} une astre ou {\`a} un
amas d'astres. 

\begin{example} \label{pdpi}
Dans la Figure \ref{fig:Suitecl} est repr{\'e}sent{\'e} un exemple sch{\'e}matique
de suite d'{\'e}clatements associ{\'e}e {\`a} une constellation $\mathcal{C}$ de 15 astres
infiniment voisins de $O$, y compris $O$. Les fl{\`e}ches indiquent des morphismes
d'{\'e}clatements successifs, $\pi_i$ {\'e}tant l'{\'e}clatement 
simultan{\'e} des astres de hauteur $i$. Les astres {\'e}clat{\'e}s sont
indiqu{\'e}s par des 
ast{\'e}risques. Les composantes du diviseur exceptionnel apparues {\`a}
chaque {\'e}tape sont indiqu{\'e}es en traits gras. On num{\'e}rote les
astres diff{\'e}rents de $O$ par $A_1,...,A_{14}$. La composante
cr{\'e}{\'e}e par l'{\'e}clatement de $A_i$ est not{\'e}e $E_i$, la m{\^e}me notation
servant pour toutes
ses transform{\'e}es strictes. L'application $p_D: \mathcal{C} \rightarrow
\mathcal{C}$ est donn{\'e}e par : 

$$ \begin{array}{ccccccccccccccc}
     O & A_1 & A_2 & A_3 & A_4 & A_5 & A_6 & A_7 & A_8 & A_9 & A_{10}
     & A_{11} & A_{12} & A_{13} & A_{14}\\
     \downarrow & \downarrow &  \downarrow & \downarrow & \downarrow &
     \downarrow &  \downarrow & \downarrow &  \downarrow & \downarrow
     &  \downarrow & \downarrow &  \downarrow & \downarrow &
     \downarrow \\
     O & O & O & A_1 & A_1 & A_1 & A_2 & A_2 & A_4 & A_3 & A_3 & A_6 &
     A_7 & A_7  & A_6
   \end{array}  $$

Les astres libres sont $A_1, A_2, A_4, A_5, A_7, A_{12}, A_{13},
A_{14}$ et les satellites sont $A_3, A_6, A_8, A_9, A_{10},
A_{11}$. Pour ces derniers, l'application $p_I$ est donn{\'e}e par :
  $$ \begin{array}{cccccc}
     A_3 & A_6 & A_8 & A_9 & A_{10} & A_{11} \\
     \downarrow & \downarrow &  \downarrow & \downarrow & \downarrow &
     \downarrow  \\
     O & O & A_1 & A_1 & O & A_2 
   \end{array}  $$

\end{example}

\section{Le diagramme d'Enriques 
  et le graphe dual d'une constellation}

Dans cette section je rapp\`ele les deux principaux codages de la
combinatoire d'une constellation : son \emph{diagramme d'Enriques} et son
\emph{graphe dual}. 
\medskip

Voici d'abord la d{\'e}finition du diagramme d'Enriques :

\begin{definition} \label{enriq}
  Soit $\mathcal{C}$ une constellation finie. Son {\bf diagramme
    d'Enriques} $\mathcal{E}(\mathcal{C})$  est le graphe d{\'e}cor{\'e}
  enracin{\'e} d{\'e}fini de la mani{\`e}re  suivante :

  $\bullet$ ses sommets sont en bijection avec les astres de
  $\mathcal{C}$ ; sa racine correspond {\`a} $O$. 

  $\bullet$ deux sommets sont reli{\'e}s par une ar{\^e}te lorsqu'ils
  repr{\'e}sentent des astres dont l'un est directement voisin de
  l'autre. 

  $\bullet$ une ar{\^e}te est \textbf{courbe} si l'astre sup{\'e}rieur est
  libre ; sinon elle est \textbf{droite} ; deux ar{\^e}tes droites successives
  vont tout droit lorsque leurs deux astres sup{\'e}rieurs ont
  m{\^e}me pr{\'e}d{\'e}cesseur indirect ; sinon, elles forment une ligne
  bris{\'e}e ; une ar{\^e}te droite sortant d'une ar{\^e}te courbe a la m{\^e}me
  tangente que celle-ci au sommet commun ; tous les autres couples
  d'ar{\^e}tes successives forment une ligne bris{\'e}e.
\end{definition}

Les r{\`e}gles pr{\'e}c{\'e}dentes ont {\'e}t{\'e} d{\'e}crites dans \cite{EC 18}
afin de permettre de dessiner le diagramme dans le plan. Mais il faut
bien voir qu'elles d{\'e}crivent en fait uniquement une structure
suppl{\'e}mentaire sur un graphe abstrait, sans privil{\'e}gier un 
plongement plan par rapport {\`a} un autre. 

La distance g{\'e}od{\'e}sique d'un sommet de
$\mathcal{E}(\mathcal{C})$ {\`a} 
$O$ est {\'e}gale {\`a} la hauteur de l'astre correspondant. Orientons
chaque ar{\^e}te de son sommet le plus bas (dit \emph{sommet initial})
vers son sommet le plus haut (dit \emph{sommet terminal}). Il est
imm{\'e}diat de voir que l'on a les
r{\`e}gles suivantes pour lire sur le diagramme d'Enriques
$\mathcal{E}(\mathcal{C})$ les fonctions $p_D$ et $p_I$ :

\begin{proposition} \label{lecture}
  Un sommet $A$ de $\mathcal{E}(\mathcal{C})$ correspond {\`a}
      un astre satellite si et seulement si l'ar{\^e}te qui y aboutit
      est droite. Dans ce cas, $p_D(A)=B$, o{\`u} $B$ est le sommet
      initial de l'ar{\^e}te aboutissant {\`a} $A$ et $p_I(A)=C$, o{\`u} :
  \begin{enumerate}
    \item $C$ est le sommet initial de l'ar{\^e}te aboutissant {\`a} $B$
      si cette ar{\^e}te et $BA$ forment une ligne bris{\'e}e ;

    \item sinon,  $C$ est le sommet le plus bas sur la g{\'e}od{\'e}sique
      joignant $B$ {\`a} $O$, tel que la g{\'e}od{\'e}sique $CB$ ne soit
      pas bris{\'e}e.
  \end{enumerate}
\end{proposition}

\begin{figure}[h!] 
\vspace*{6mm}
\labellist \small\hair 2pt 
  \pinlabel{$O$} at 280 -10
  \pinlabel{$A_1$} at 214 97
  \pinlabel{$A_2$} at 346 102
  \pinlabel{$A_3$} at 160 217
  \pinlabel{$A_4$} at 83 93
  \pinlabel{$A_5$} at 266 152
  \pinlabel{$A_6$} at 395 192
  \pinlabel{$A_7$} at 463 110
  \pinlabel{$A_8$} at 12 49
  \pinlabel{$A_9$} at 241 262
  \pinlabel{$A_{10}$} at 79 293
  \pinlabel{$A_{11}$} at 316 249
  \pinlabel{$A_{12}$} at 482 230
  \pinlabel{$A_{13}$} at 545 100
  \pinlabel{$A_{14}$} at 495 292
\endlabellist 
\centering 
\includegraphics[scale=0.40]{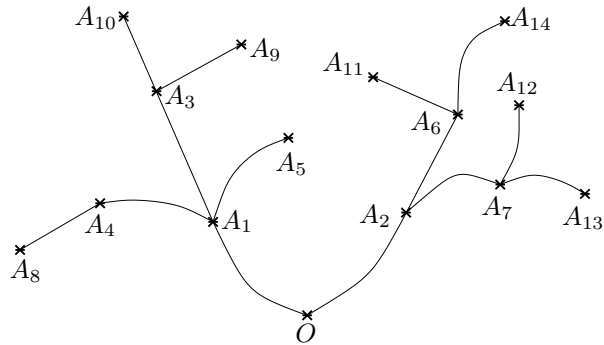} 
\vspace*{10mm} \caption{Le diagramme d'Enriques de la constellation
  de la Figure  \ref{fig:Suitecl}}  
\label{fig:Enriqex} 
\end{figure}

\begin{example}
Dans la figure \ref{fig:Enriqex} est repr{\'e}sent{\'e} le
diagramme d'Enriques de la constellation de la Figure
\ref{fig:Suitecl}. Gr{\^a}ce {\`a} la Proposition \ref{lecture}, on
v{\'e}rifie les valeurs de $p_I$ donn{\'e}es dans l'Exemple \ref{pdpi}. 

\end{example}

Un deuxi{\`e}me type de graphe d{\'e}cor{\'e} est utilis{\'e} pour
repr{\'e}senter la combinatoire d'une constellation. Sa d{\'e}finition utilise
le morphisme (\ref{eclconst}) : 

\begin{definition} \label{dualgraph}
   Soit $\mathcal{C}$ une constellation finie. Son {\bf graphe dual}
   $\mathcal{D}(\mathcal{C})$ est le graphe dual du diviseur r{\'e}duit
   $E_{\mathcal{C}}$ de la 
   surface lisse $S_{\mathcal{C}}$, chaque sommet {\'e}tant pond{\'e}r{\'e}
   par l'auto-intersection de la composante irr{\'e}ductible qui lui
   correspond. 
\end{definition}

Rappelons que ces auto-intersections peuvent se calculer
r{\'e}cursivement lors du processus d'{\'e}clatements, en utilisant le fait
que lorsqu'on {\'e}clate un point sur une courbe compacte lisse $\Gamma$, 
et que l'on d{\'e}signe par $\Gamma_1$ sa transform{\'e}e stricte, on a :
  $$\Gamma_1^2 = \Gamma^2 -1. $$

\begin{figure}[h!] 
\vspace*{16mm}
\labellist \small\hair 2pt 
\pinlabel{$E_4$} at 8 -16
\pinlabel{$E_8$} at 80 20
\pinlabel{$E_1$} at 152 60
\pinlabel{$E_9$} at 193 60
\pinlabel{$E_3$} at 240 60
\pinlabel{$E_{10}$} at 285 60
\pinlabel{$E_0$} at 324 60
\pinlabel{$E_6$} at 355 60
\pinlabel{$E_{11}$} at 409 60
\pinlabel{$E_2$} at 454 60
\pinlabel{$E_7$} at 494 60
\pinlabel{$E_{12}$} at 543 21
\pinlabel{$E_{13}$} at 543 137
\pinlabel{$E_5$} at 77 137
\pinlabel{$E_{14}$} at 370 5

\pinlabel{$-2$} at 8 24
\pinlabel{$-1$} at 80 60
\pinlabel{$-6$} at 152 100
\pinlabel{$-1$} at 193 100
\pinlabel{$-3$} at 240 100
\pinlabel{$-1$} at 285 100
\pinlabel{$-6$} at 324 100
\pinlabel{$-3$} at 369 100
\pinlabel{$-1$} at 409 100
\pinlabel{$-4$} at 454 100
\pinlabel{$-3$} at 494 100
\pinlabel{$-1$} at 543 61
\pinlabel{$-1$} at 543 97
\pinlabel{$-1$} at 77 97
\pinlabel{$-1$} at 390 25
\endlabellist 
\centering 
\includegraphics[scale=0.60]{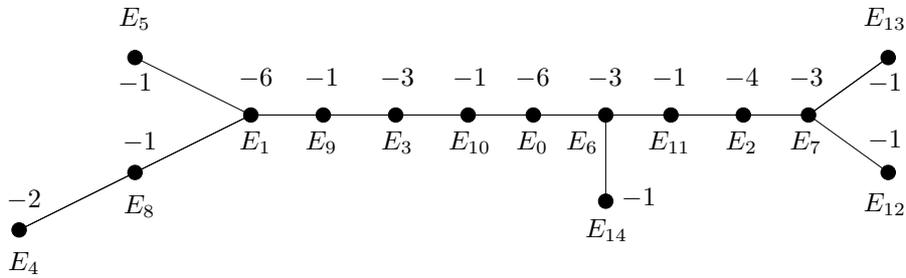} 
\vspace*{2mm} \caption{Le graphe dual de la constellation de la Figure
  \ref{fig:Suitecl}} 
\label{fig:Dualex} 
\end{figure}

\begin{example}
Dans la Figure \ref{fig:Dualex} est repr{\'e}sent{\'e} le graphe
dual associ{\'e} {\`a} la constellation de la Figure \ref{fig:Suitecl}. Les
sommets sont 
num{\'e}rot{\'e}s par les composantes irr{\'e}ductibles de $E_{\mathcal{C}}$ qu'ils
repr{\'e}sentent, et ils sont pond{\'e}r{\'e}s par les auto-intersections
respectives. 
\end{example}

\begin{remark}
   Le graphe dual d'un diviseur r{\'e}duit sur une surface lisse est consid{\'e}r{\'e} en
   passant par Du Val \cite{DV 
     36}. Mais il ne semble avoir commenc{\'e} {\`a} {\^e}tre utilis{\'e}
   syst{\'e}matiquement qu'{\`a} la suite de l'article \cite{M 61} de
   Mumford et de la 
   pr{\'e}sentation \cite{H 63} qui en a {\'e}t{\'e} faite par Hirzebruch. 
\end{remark}

Les sommets des graphes $\mathcal{E}(\mathcal{C})$ et
$\mathcal{D}(\mathcal{C})$  sont en bijection naturelle : on associe
au sommet de $\mathcal{E}(\mathcal{C})$ repr{\'e}sentant l'astre $A$
le sommet de $\mathcal{D}(\mathcal{C})$ repr{\'e}sentant la courbe
exceptionnelle obtenue en {\'e}clatant $A$. Cette bijection ne respecte
pas les structures des deux graphes. En fait, ceux-ci ne sont en
g{\'e}n{\'e}ral m{\^e}me
pas abstraitement isomorphes, comme on le voit en comparant les
figures \ref{fig:Enriqex}  et \ref{fig:Dualex}. N{\'e}anmoins, ils
codent tous les deux la m{\^e}me information : il est possible de donner
des \emph{algorithmes} permettant de passer de l'un {\`a} l'autre (voir
\cite{CA 00} ou \cite{CC 05}). 

L'un des buts de cet article est de faciliter la
compr{\'e}hension \emph{g{\'e}om{\'e}trique} de la relation entre les deux
graphes. 

\medskip
L'id{\'e}e de base est de repr{\'e}senter chaque astre d'une constellation
par \emph{deux} points distincts : un premier le repr{\'e}sentant en
tant que point ferm{\'e} sur l'un des mod{\`e}les, et un deuxi{\`e}me  
repr{\'e}sentant le diviseur
exceptionnel cr{\'e}{\'e} par l'{\'e}clatement de ce point ferm{\'e}. De
plus, chaque fois que le point ferm{\'e} sera vu comme intersection de
deux courbes lisses transverses, on aura un triangle affine
canoniquement associ{\'e} {\`a} ce diviseur {\`a} croisements normaux et un
plongement canonique des deux points dans le triangle. 

\begin{figure}[h!] 
\vspace*{6mm}
\labellist \small\hair 2pt 
\pinlabel{E} at 150 65
\pinlabel{E'} at 38 142
\pinlabel{A} at 52 56
\pinlabel{E} at 430 2
\pinlabel{E'} at 285 2
\pinlabel{$A$} at 356 2
\pinlabel{$E_A$} at 355 148
\pinlabel{$\Sigma$} at 158 23
\endlabellist 
\centering 
\includegraphics[scale=0.60]{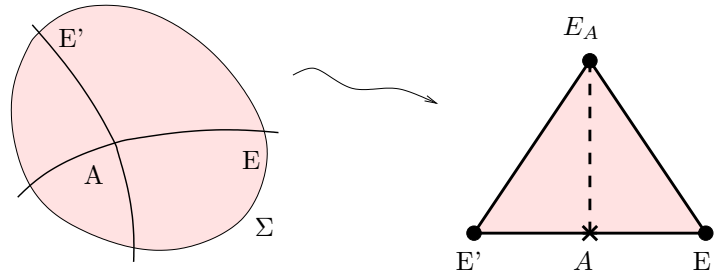} 
\vspace*{5mm} \caption{Le triangle associ{\'e} {\`a} un diviseur {\`a}
  croisements normaux} 
\label{fig:Trifond} 
\end{figure}

Expliquons cela avec plus de d{\'e}tails. 
Soit $(E \cup E', A) \hookrightarrow \Sigma$ un germe de diviseur {\`a}
croisements normaux sur une surface lisse $\Sigma$. C'est-{\`a}-dire que
$E$ et $E'$ sont deux germes en $A$ de courbes lisses
transverses. Notons par $E_A$ le diviseur exceptionnel de
l'{\'e}clatement de $A$ dans $\Sigma$. On associe au diviseur {\`a} croisements normaux
$(E \cup E', A)$ un triangle affine dont les sommets correspondent
bijectivement aux courbes $E, E', E_A$ et le milieu  du segment $[E, E']$
au point $A$. {\`A} l'astre $A$ correspondent de cette mani{\`e}re
\emph{deux} points privil{\'e}gi{\'e}s dans le triangle, $A$ et
$E_A$ (voir la Figure \ref{fig:Trifond}). On peut penser que $A$
repr{\'e}sente la courbe $E_A$ sous forme embryonnaire, et que le segment
qui les relie dans le triangle repr{\'e}sente l'embryog{\'e}n{\`e}se.

Lorsque l'on consid{\`e}re le processus d'{\'e}clatements associ{\'e} {\`a}
une constellation, on se retrouve avec une suite de germes de
diviseurs {\`a} croisements normaux : les germes des diviseurs
exceptionnels des compos{\'e}s d'{\'e}clatements aux astres satellites de
la constellation. On associe {\`a} chacun d'entre eux
un triangle comme pr{\'e}c{\'e}demment. Ces divers triangles se recollent
ensuite naturellement. Plus de soin doit {\^e}tre accord{\'e} aux astres
libres, pour lesquels on a seulement des \emph{demi-triangles}.

\medskip
\section{Construction des voilures et des cerf-volants}
  \label{constvoil}

Dans cette section j'explique les d{\'e}finitions de deux types 
de complexes simpliciaux g{\'e}om{\'e}triques bidimensionnels connexes,
les \emph{voilures} et les \emph{cerf-volants}. Ces d{\'e}finitions sont
r{\'e}cursives, par rajouts de pi{\`e}ces {\'e}l{\'e}mentaires triangulaires appel{\'e}es les
\emph{demi-voiles} et les \emph{voiles simples} et de segments
appel{\'e}s les \emph{cordes},
le tout rattach{\'e} {\`a} un segment initial appel{\'e} l'\emph{axe} du
cerf-volant.  Les structures affines des voiles simples et des
demi-voiles se recollent canoniquement, comme  expliqu{\'e} dans
la Section \ref{modaff}. Le lien avec les constellations est quant {\`a}
lui expliqu{\'e} dans la Section \ref{sousarbres}.
\medskip

La d{\'e}finition suivante introduit des termes permettant de parler en
termes intuitifs des pi{\`e}ces du jeu de construction de cerf-volants
et de leurs r{\`e}gles de recollement.

\begin{figure}[h!] 
\vspace*{6mm}
\labellist \small\hair 2pt 
\pinlabel{sommet} at 0 -10
\pinlabel{{\'e}toil{\'e}} at 0 -30
\pinlabel{sommet} at 80 -10
\pinlabel{de base} at 80 -30
\pinlabel{sommet} at 0 150
\pinlabel{terminal} at 0 130
\pinlabel{c{\^o}t{\'e}} at 80 70
\pinlabel{lat{\'e}ral} at 80 50

\pinlabel{premier} at 222 -10
\pinlabel{sommet} at 222 -30
\pinlabel{de base} at 222 -50
\pinlabel{deuxi{\`e}me} at 365 -10
\pinlabel{sommet} at 365 -30
\pinlabel{de base} at 365 -50
\pinlabel{point} at 294 -10
\pinlabel{{\'e}toil{\'e}} at 294 -30
\pinlabel{sommet} at 294 150
\pinlabel{terminal} at 294 130
\pinlabel{c{\^o}t{\'e}} at 222 70
\pinlabel{lat{\'e}ral} at 222 50
\pinlabel{c{\^o}t{\'e}} at 365 70
\pinlabel{lat{\'e}ral} at 365 50
\pinlabel{{\bf demi-voile}} at 40 -80
\pinlabel{{\bf voile simple}} at 294 -80
\endlabellist 
\centering 
\includegraphics[scale=0.60]{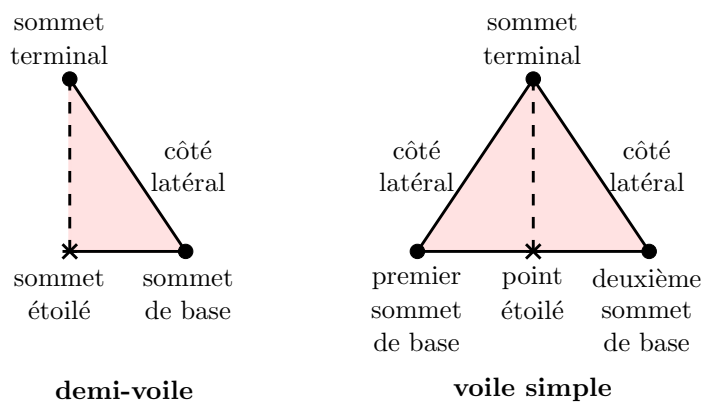} 
\vspace*{15mm} \caption{Les voiles {\'e}l{\'e}mentaires} 
\label{fig:Voilelem} 
\end{figure}

\begin{definition} (voir la Figure \ref{fig:Voilelem})
 Une {\bf demi-voile} est un triangle affine ayant un sommet
 {\bf {\'e}toil{\'e}} et deux sommets {\bf pleins}, l'un d'entre eux {\'e}tant {\bf
   de base} et l'autre {\bf terminal}.  Le {\bf c{\^o}t{\'e} lat{\'e}ral} de la
  demi-voile est celui qui joint les deux sommets pleins.

 Une {\bf voile simple} est un triangle affine dont tous les sommets
 sont {\bf pleins}, l'un d'entre eux {\'e}tant {\bf terminal} et les deux autres
 {\'e}tant {\bf de base}. De plus, ces derniers sont ordonn{\'e}s : on
 parlera du premier et du deuxi{\`e}me sommet de base. Le c{\^o}t{\'e}
 joignant les deux sommets de base est {\bf la base} et
 son milieu est {\bf le point {\'e}toil{\'e}} de la voile. Les c{\^o}t{\'e}s
 joignant le sommet  
 terminal {\`a} l'un des sommets de base sont appel{\'e}s {\bf c{\^o}t{\'e}s
   lat{\'e}raux}. 

On dira que les demi-voiles et les voiles simples sont les {\bf voiles
  {\'e}l{\'e}mentaires}. 
\end{definition}

Le vocabulaire pr{\'e}c{\'e}dent est motiv{\'e} par le fait qu'une
demi-voile est affinement isomorphe aux moiti{\'e}s des voiles
simples obtenues en joignant le point {\'e}toil{\'e} au sommet terminal
par le segment affine qui les relie. Dans les dessins, nous
repr{\'e}senterons ces segments par des traits hachur{\'e}s. Ce sont eux
qui mat{\'e}rialiseront la correspondance bijective naturelle entre
points {\'e}toil{\'e}s et sommets pleins dans les complexes simpliciaux
que nous construirons. 

Les cordes et l'axe ont 
aussi leurs sommets d{\'e}cor{\'e}s en types :

\begin{definition}
 Une {\bf corde} est un segment affine dont l'un des sommets est  {\bf
   initial} et l'autre {\bf final}. L'{\bf axe} est un segment affine
 dont l'un des sommets est {\bf {\'e}toil{\'e}} et l'autre est {\bf
   plein}.  
\end{definition}

\begin{figure}[h!] 
\vspace*{6mm}
\labellist \small\hair 2pt 
\pinlabel{recollement} at 150 115
\pinlabel{des sommets} at 150 95
\pinlabel{voilure} at 335 90
\pinlabel{pr{\'e}existante} at 335 70
\endlabellist 
\centering 
\includegraphics[scale=0.50]{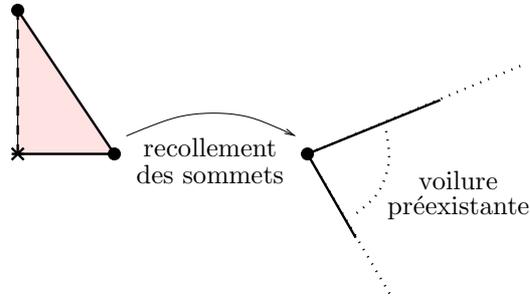} 
\caption{Rattachement d'une demi-voile} 
\label{fig:Ratdemi} 
\end{figure}

\begin{figure}[h!] 
\vspace*{10mm}
\labellist \small\hair 2pt 
\pinlabel{recollement} at 265 106
\pinlabel{voile {\'e}l{\'e}mentaire} at 113 38
\pinlabel{de rattachement} at 113 18
\pinlabel{voilure} at 25 120
\pinlabel{pr{\'e}existante} at 25 100

\pinlabel{deuxi{\`e}me} at 60 250
\pinlabel{sommet de base} at 60 230
\pinlabel{sommet terminal} at -35 180
\pinlabel{de la voile de rattachement} at -35 160

\pinlabel{disparition} at 292 68
\pinlabel{de ce c{\^o}t{\'e}} at 292 48
\pinlabel{lat{\'e}ral} at 292 28
\pinlabel{nouveaux} at 300 340
\pinlabel{c{\^o}t{\'e}s} at 300 320
\pinlabel{lat{\'e}raux} at 300 300
\endlabellist 
\centering 
\includegraphics[scale=0.50]{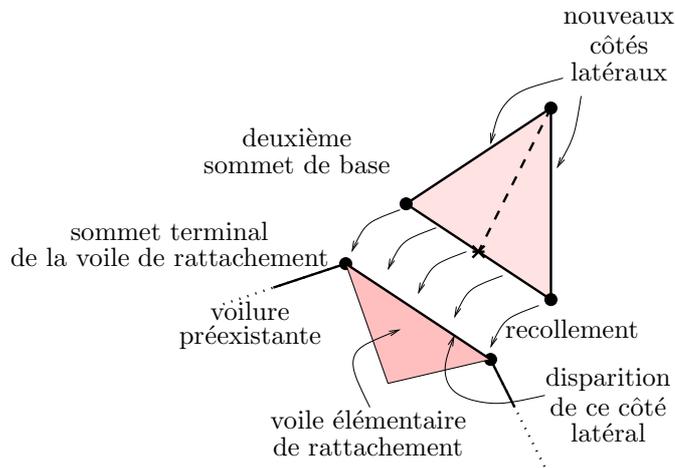} 
\caption{Rattachement d'une voile simple} 
\label{fig:Ratsimple} 
\end{figure}

En ayant {\`a} notre disposition un kit de construction form{\'e} d'un seul axe,
mais de demi-voiles et de voiles simples {\`a} volont{\'e}, nous pouvons 
assembler des \emph{voilures} plus compliqu{\'e}es par un processus de
construction  dont les {\'e}tapes {\'e}l{\'e}mentaires sont :

\begin{enumerate}
  \item On part de l'axe, consid{\'e}r{\'e} comme une voilure d{\'e}g{\'e}n{\'e}r{\'e}e.
  
  \item Si $\mathcal{V}$ est une voilure d{\'e}j{\`a} construite, on peut au choix :
  
    \begin{enumerate}
         \item prendre une nouvelle demi-voile et recoller son sommet de
           base {\`a}  un sommet plein de $\mathcal{V}$ (voir la Figure
           \ref{fig:Ratdemi});  
     
        \item prendre une nouvelle voile simple et recoller sa base
          {\`a} un c{\^o}t{\'e} lat{\'e}ral de $\mathcal{V}$  par l'unique isomorphisme
                  affine qui envoie le deuxi{\`e}me sommet de base de la
                  nouvelle voile  simple sur le sommet  
                  terminal de la voile {\'e}l{\'e}mentaire {\`a} laquelle
                  appartient le c{\^o}t{\'e} lat{\'e}ral. On remplace ainsi un c{\^o}t{\'e}
                  lat{\'e}ral de l'ancienne voilure par deux nouveaux c{\^o}t{\'e}s
                  lat{\'e}raux    
                  (voir la Figure \ref{fig:Ratsimple}).
    \end{enumerate}
\end{enumerate}

Si l'on enl{\`e}ve les sommets d'une voilure $\mathcal{V}$, celle-ci se
d{\'e}compose en 
composantes connexes. Appelons \emph{voile compl{\`e}te} l'adh{\'e}rence de chacune
de ces composantes connexes dans la voilure $\mathcal{V}$. 
L'une de
ces voiles compl{\`e}tes est l'axe de la 
voilure. Chacune des autres voiles compl{\`e}tes est un complexe
simplicial purement bidimensionnel hom{\'e}omorphe {\`a} un disque, ayant 
une unique demi-voile, le reste des triangles {\'e}tant des voiles
simples. 

\begin{example}
Dans la Figure \ref{fig:Voilex} est repr{\'e}sent{\'e}e une
voilure. Cette voilure a 8 voiles compl{\`e}tes en dehors de l'axe, 5
d'entre elles  {\'e}tant r{\'e}duites {\`a} une demi-voile, les 3 restantes 
ayant 1, 2 et respectivement 3 voiles simples dans leur constitution. 
\end{example}

\medskip

Chaque voile compl{\`e}te \emph{s'oriente canoniquement} de la mani{\`e}re
suivante :

$\bullet$ on oriente la demi-voile initiale en choisissant l'ordre
suivant des sommets : sommet de base, sommet terminal, sommet {\'e}toil{\'e} ;

$\bullet$ on propage cette orientation par continuit{\'e} {\`a} toute la
voile compl{\`e}te. 

Ceci permet de parler de \emph{c{\^o}t{\'e} lat{\'e}ral droit} et de
\emph{c{\^o}t{\'e} lat{\'e}ral gauche} de chaque voile simple : le c{\^o}t{\'e}
droit est celui que l'on rencontre en tournant positivement lorsque
l'on sort de la base. Lorsque l'on recolle une voile simple {\`a} une
autre, on peut donc dire si le recollement se fait \emph{sur le
  c{\^o}t{\'e} droit} ou \emph{sur le c{\^o}t{\'e} gauche}. Consid{\'e}rons une suite
de voiles simples $(\tau_1,..., \tau_n)$ recoll{\'e}es les unes aux autres
dans cet ordre. Si les recollements 
sont toujours effectu{\'e}s du m{\^e}me c{\^o}t{\'e}, on dira que ces voiles
\emph{tournent dans le m{\^e}me sens}. Cette notion peut s'{\'e}tendre au
cas o{\`u} l'on part d'une demi-voile et que le recollement se fait
continuellement sur le c{\^o}t{\'e} droit. 
\medskip

\begin{figure}[h!] 
\vspace*{6mm}
\labellist \small\hair 2pt 
\pinlabel{recollement} at 150 155
\pinlabel{voilure} at 355 110
\pinlabel{pr{\'e}existante} at 355 90
\pinlabel{{\bf corde libre}} at 160 -10
\pinlabel{sommet initial} at 355 30
\pinlabel{sommet final} at -15 60
\endlabellist 
\centering 
\includegraphics[scale=0.50]{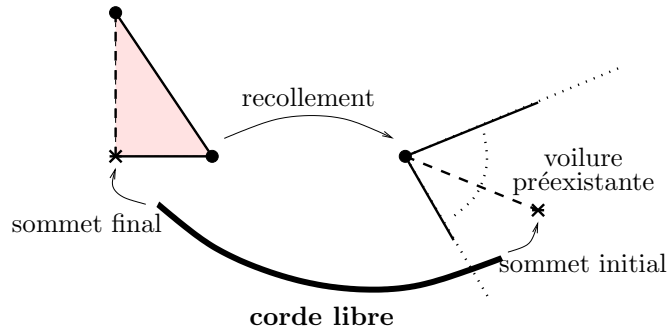} 
\vspace*{2mm} \caption{Rattachement d'une corde libre} 
\label{fig:Ratcorde1} 
\end{figure}

{\`A} chaque voilure $\mathcal{V}$ on associe canoniquement un 
\emph{cerf-volant} $\mathcal{KV}$. Pour cela, lors de la construction
de $\mathcal{V}$  \emph{on attache une 
corde en m{\^e}me temps qu'on recolle une voile {\'e}l{\'e}mentaire}. Ceci se fait
de la mani{\`e}re suivante : 

(a) si l'on recolle une \emph{demi-voile}, on attache aussi une corde en
identifiant son sommet final au sommet {\'e}toil{\'e} de la demi-voile et
son sommet initial au sommet {\'e}toil{\'e} ou au point {\'e}toil{\'e} qui
correspond au sommet plein auquel a {\'e}t{\'e} attach{\'e}e la demi-voile.  On
dira qu'il s'agit d'une {\bf corde  libre} du cerf-volant (voir
la Figure \ref{fig:Ratcorde1}). 

\begin{figure}[h!] 
\vspace*{6mm}
\labellist \small\hair 2pt 
\pinlabel{M} at 143 44
\pinlabel{{\bf corde satellite}} at 50 89
\pinlabel{C} at 220 82
\pinlabel{B} at 80 165
\endlabellist 
\centering 
\includegraphics[scale=0.50]{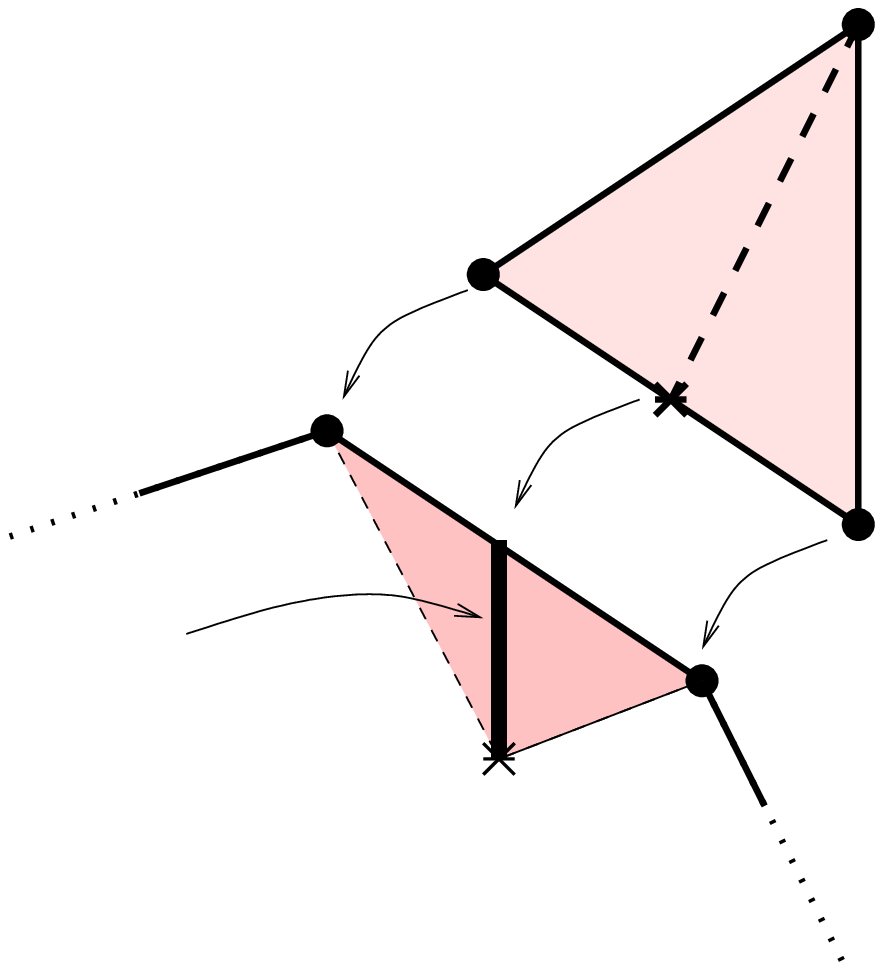} 
\caption{Rattachement d'une corde satellite} 
\label{fig:Ratcorde2} 
\end{figure}

(b) si l'on recolle une \emph{voile simple}, notons par $\tau$ la voile
{\'e}l{\'e}mentaire de $\mathcal{V}$ {\`a} laquelle on la rattache, par $B$
son sommet terminal, par $C$ le sommet de base tel que le segment $BC$
soit celui le long duquel la nouvelle voile simple est attach{\'e}e et
par $M$ le point {\'e}toil{\'e} de $\tau$. On recolle alors une corde le
long du segment joignant $M$ au milieu du segment $BC$, par un
isomorphisme affine qui envoie le sommet initial de la corde sur $M$. 
On dira qu'il s'agit d'une {\bf corde  satellite} du cerf-volant 
(voir la Figure \ref{fig:Ratcorde2}).  

\begin{example}  Dans la Figure \ref{fig:Cerfconsex} est
  repr{\'e}sent{\'e} le cerf-volant associ{\'e} {\`a} la voilure de la Figure
  \ref{fig:Voilex}. Comme la figure est plane, certaines cordes libres
  - que l'on repr{\'e}sentera toujours par des lignes courbes -
  sont parfois oblig{\'e}es pour des raisons topologiques d'{\^e}tre
  repr{\'e}sent{\'e}es intersectant la voilure ailleurs qu'en leurs
  extr{\'e}mit{\'e}s. C'est ici le cas de la corde joignant $A_6$ et
  $A_{14}$. Quant aux cordes $A_1 A_5$ et $A_7 A_{12}$, on aurait pu
  les dessiner sans de telles intersections, mais on a pr{\'e}f{\'e}r{\'e}
  montrer que ces intersections suppl{\'e}mentaires ne nuisent pas
  tellement {\`a} la lisibilit{\'e} de la figure, une fois l'\oe il
  entra{\^\i}n{\'e}. 
\end{example}

\medskip 

\section{Mod{\`e}le affine canonique des voiles 
    compl{\`e}tes } \label{modaff}

Dans cette section j'explique comment associer canoniquement {\`a}
chaque base d'un r{\'e}seau bidimensionnel un
complexe simplicial de dimension deux plong{\'e} dans le c{\^o}ne convexe
engendr{\'e} - son \emph{lotus}. Puis j'explique comment plonger
canoniquement chaque voile compl{\`e}te d'une voilure dans le lotus. Ce
plongement 
d{\'e}finit un \emph{recollement canonique} des structures affines des voiles
{\'e}l{\'e}mentaires constituant chaque voile compl{\`e}te. 
\medskip

Consid{\'e}rons un r{\'e}seau bidimensionnel $N$ (c'est-{\`a}-dire un groupe
ab{\'e}lien libre de rang $2$) et une base $(e_1, e_2)$ de ce
r{\'e}seau. Notons par $\sigma(e_1, e_2)$ le c{\^o}ne convexe engendr{\'e}
par la base dans 
l'espace vectoriel r{\'e}el $N_{\R}:=N \otimes_{\Z} \R$
associ{\'e}. Notons par $\tau(e_1, e_2)$ le triangle qui est contenu dans le
plan r{\'e}el $N_{\R}$ et qui joint les points $e_1, e_2, e_1 + e_2$. 

\begin{figure}[h!] 
\vspace*{5mm}
\labellist \small\hair 2pt 
  \pinlabel{0} at -5 20
  \pinlabel{$e_1$} at 21 20
  \pinlabel{$e_2$} at -5 40
  \pinlabel{$\tau(e_1, e_2)$} at 25 -5
\endlabellist 
\centering 
\includegraphics[scale=1.50]{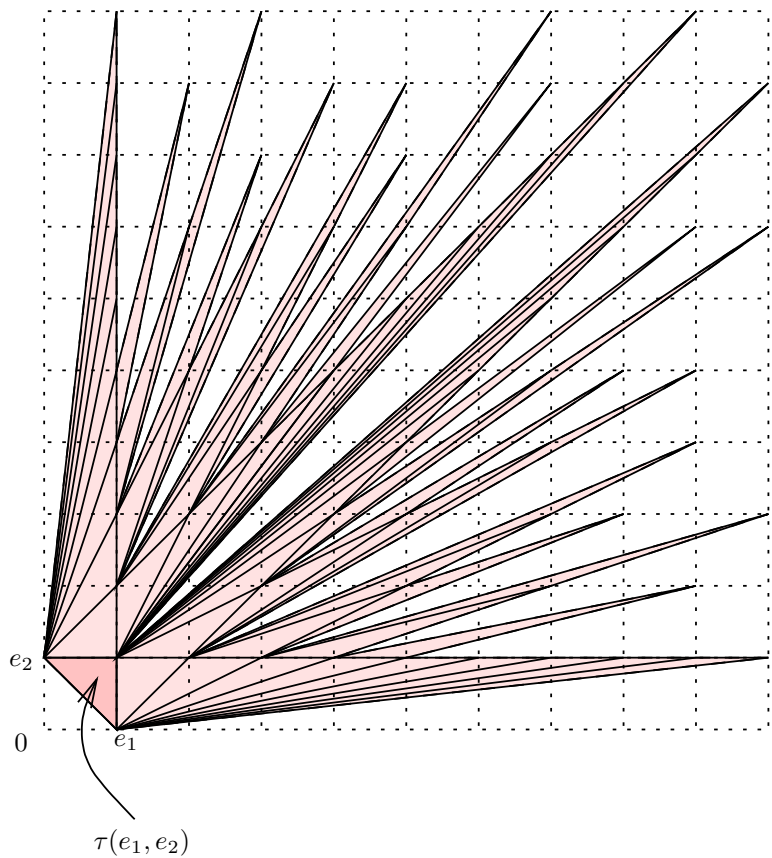} 
\vspace*{5mm} \caption{Le lotus $\mathcal{L}(e_1, e_2)$} 
\label{fig:Lotus} 
\end{figure}

Cette construction peut {\^e}tre ensuite r{\'e}p{\'e}t{\'e}e {\`a} partir des
bases $(e_1, e_1 +e_2)$ et $(e_2, e_1 +e_2)$ de $N$. Ainsi, de proche
en proche, on construit un complexe simplicial infini plong{\'e} dans le
c{\^o}ne $\sigma(e_1, e_2)$ : {\`a} la $n$-{\`e}me {\'e}tape de construction, on
rajoute $2^n$ triangles {\`a} ceux d{\'e}j{\`a} construits. Dans la Figure
\ref{fig:Lotus} est  
repr{\'e}sent{\'e}e l'union de tous les triangles de ce complexe simplicial
contenus dans le parall{\'e}logramme engendr{\'e} par $10 e_1, 10 e_2$. {\`A}
cause de cette forme, Bernard Teissier m'a sugg{\'e}r{\'e} :

\begin{definition}
 Le complexe simplicial pr{\'e}c{\'e}dent est appel{\'e} {\bf le lotus}
 $\mathcal{L}(e_1, e_2)$ associ{\'e} au 
 c{\^o}ne $\sigma(e_1, e_2)$. 
\end{definition}

Je dirai aussi, en filant la
m{\'e}taphore, que les triangles sont les \emph{p{\'e}tales} du lotus. 

Bien s{\^u}r, {\`a}
transformations affines pr{\'e}servant les r{\'e}seaux pr{\`e}s, il n'y a
qu'un seul lotus.

\begin{remark}
  Consid{\`e}rons la vari{\'e}t{\'e} torique affine de r{\'e}seau des
  poids $N$ et d'{\'e}ventail form{\'e} par le c{\^o}ne $\sigma(e_1, e_2)$ et ses
  faces. Elle est isomorphe {\`a} $\C^2$, munie de sa structure
  torique canonique. Consid{\'e}rons la suite des {\'e}clatements des
  orbites de dimension $0$. C'est une suite de morphismes toriques,
  obtenue en subdivisant successivement les c{\^o}nes de dimension $2$
  de l'{\'e}ventail de l'{\'e}tape pr{\'e}c{\'e}dente. Si on dessine {\`a}
  chaque fois le triangle ayant comme sommets les deux vecteurs
  primitifs des ar{\^e}tes d'un c{\^o}ne de dimension $2$ et celui de l'ar{\^e}te le
  subdivisant, on obtient exactement les p{\'e}tales du lotus. Le lotus
  permet donc 
  de visualiser d'un seul coup d'\oe il la suite infinie des morphismes
  d'{\'e}clatement des orbites de dimension $0$. 
\end{remark}

Consid{\'e}rons {\`a} pr{\'e}sent une voile compl{\`e}te d'une
voilure, diff{\'e}rente de l'axe. Elle peut se plonger canoniquement dans
le lotus $\mathcal{L}(e_1, e_2)$ par une application qui identifie
chaque voile simple {\`a} un p{\'e}tale : 

$\bullet$ on plonge la demi-voile dans $\sigma(e_1, e_2)$ par l'unique
transformation affine qui envoie le sommet {\'e}toil{\'e} en
$\frac{1}{2}(e_1 + e_2)$, le sommet de base en $e_1$ et le sommet
terminal en $e_1 + e_2$. 

$\bullet$ on plonge la voile simple recoll{\'e}e le long du c{\^o}t{\'e}
lat{\'e}ral de la demi-voile en envoyant son premier sommet de base sur
$e_1$, son deuxi{\`e}me sommet de base sur $e_1 + e_2$ et son sommet
terminal sur $2 e_1 + e_2$.

$\bullet$  chaque nouvelle
voile simple se plonge en respectant les incidences.

Gr{\^a}ce {\`a} ce plongement, on obtient \emph{une structure affine enti{\`e}re
  canonique}
sur chaque voile compl{\`e}te. Le fait que la structure soit
\emph{enti{\`e}re} signifie ici simplement que l'on sait dire quels sont les
points entiers : ce sont exactement les sommets des voiles
simples. Mais le sommet {\'e}toil{\'e} de la demi-voile initiale n'est que
demi-entier !

\begin{figure}[h!] 
\vspace*{6mm}
\labellist \small\hair 2pt 
  \pinlabel{0} at 300 -5
\endlabellist 
\centering 
\includegraphics[scale=0.60]{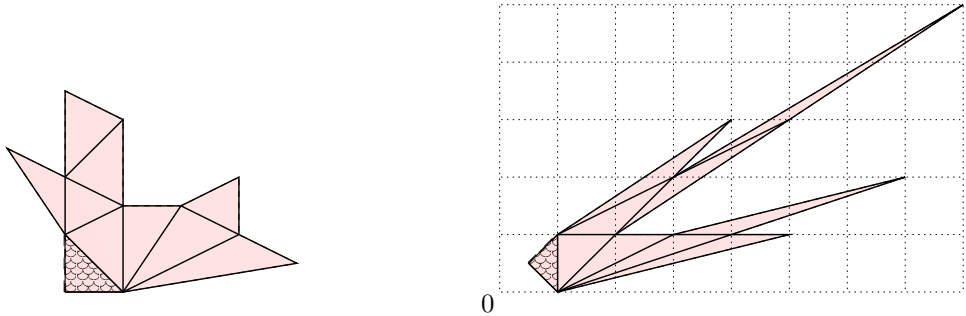} 
\vspace*{5mm} \caption{Une voile compl{\`e}te et son plongement affine canonique
  dans le lotus} 
\label{fig:Plongaf} 
\end{figure}

\begin{example}
Dans la Figure \ref{fig:Plongaf} est repr{\'e}sent{\'e} un exemple
de voile compl{\`e}te et son plongement affine canonique dans le
lotus. L'unique demi-voile est mise en {\'e}vidence {\`a} l'aide d'un
motif sp{\'e}cial.
\end{example}

Gr{\^a}ce {\`a} l'existence de ce plongement canonique, cela a un sens de
dire qu'une ligne polygonale contenue dans une voile compl{\`e}te \emph{va tout
  droit} (c'est-{\`a}-dire que c'est une g{\'e}od{\'e}sique pour la
structure affine du recollement) ou non. Ceci permet d'exprimer en
termes affines le fait qu'une suite de p{\'e}tales tourne dans le m{\^e}me
sens :

\begin{proposition} \label{toutdroit}
  {\`A} l'int{\'e}rieur de l'une des voiles compl{\`e}tes d'une voilure, on consid{\`e}re
  une suite de voiles {\'e}l{\'e}mentaires $\tau_1,..., \tau_n$ telle que
  deux voiles successives soient adjacentes et construites dans cet
  ordre. Alors elles tournent toujours dans le
  m{\^e}me sens si et seulement si la suite des cordes satellites associ{\'e}es est une
  g{\'e}od{\'e}sique pour la structure affine 
  canonique de la voile compl{\`e}te. 
\end{proposition}

\begin{figure}[h!] 
\vspace*{6mm}
\labellist \small\hair 2pt 
  \pinlabel{0} at 300 -5
\endlabellist 
\centering 
\includegraphics[scale=0.60]{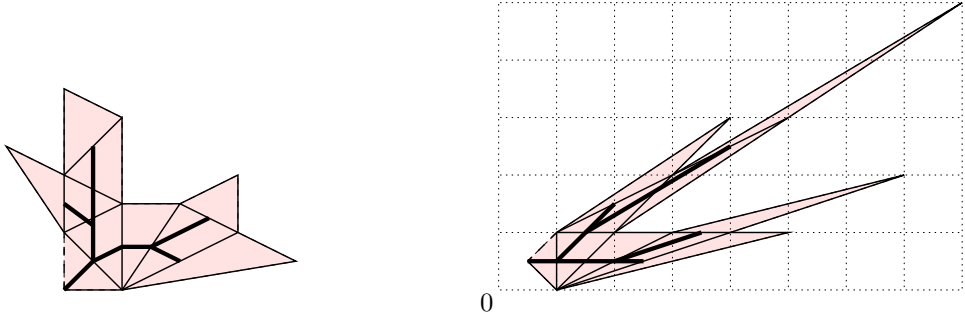} 
\vspace*{5mm} \caption{Cordes satellites d'une voile compl{\`e}te} 
\label{fig:Toudroit} 
\end{figure}

\begin{example}
Cette proposition est illustr{\'e}e dans la Figure \ref{fig:Toudroit}
pour la voile compl{\`e}te de la Figure \ref{fig:Plongaf}.  La proposition
pr{\'e}c{\'e}dente permet de 
rep{\'e}rer sur la voile compl{\`e}te repr{\'e}sent{\'e}e combinatoirement ({\`a} gauche)
les cordes satellites align{\'e}es dans le plongement canonique dans le
lotus ({\`a} droite).  Il est important de savoir faire cette
reconnaissance sur une voilure d{\'e}form{\'e}e, car lorsque le nombre de
voiles croit, tr{\`e}s rapidement le plongement affine canonique devient
impossible {\`a} dessiner {\`a} cause de l'allongement des triangles. 
\end{example}

\begin{remark}
  Consid{\'e}rons {\`a} l'int{\'e}rieur de chaque p{\'e}tale du lotus
  $\mathcal{L}(e_1, e_2)$ les deux segments joignant le 
  milieu de sa base aux milieux des c{\^o}t{\'e}s. Leur
  union est un arbre binaire plong{\'e} dans $\sigma(e_1,
  e_2)$. D{\'e}corons chaque sommet {\'e}toil{\'e} $\frac{1}{2}(a_1 e_1 +
  a_2 e_2)$ par la fraction $\frac{a_1}{a_2}$. On obtient ainsi un arbre
  isomorphe {\`a} \emph{l'arbre de Stern-Brocot} d{\'e}crit dans
  \cite[page 116]{GKP 94}. Ce dernier arbre repr{\'e}sente
  g{\'e}om{\'e}triquement la construction des suites de Farey par
  it{\'e}ration de l'op{\'e}ration $(\frac{m}{n}, \frac{m'}{n'})
  \rightarrow \frac{m+ m'}{n + n'}$, en partant de la suite $(\frac{0}{1},
  \frac{1}{0})$, et en consid{\'e}rant {\`a} chaque {\'e}tape les couples
    successifs de la suite construite {\`a} l'{\'e}tape pr{\'e}c{\'e}dente. Il est
    construit en 
    reliant chaque nouvelle fraction {\`a} celle ayant {\'e}t{\'e} cr{\'e}{\'e}e en
    dernier parmi les deux lui ayant donn{\'e} naissance. 
\end{remark}

\section{Le diagramme d'Enriques et le graphe dual 
  d'une constellation comme sous-arbres du 
cerf-volant} \label{sousarbres}

Dans cette section j'explique comment associer {\`a} chaque
constellation finie une voilure, donc aussi un cerf-volant. Puis je 
montre comment lire le graphe dual de la constellation {\`a}
partir de la voilure et le diagramme d'Enriques {\`a} partir du
cerf-volant (voir Th{\'e}or{\`e}me \ref{isos}). En fait le graphe dual est
canoniquement isomorphe {\`a} une partie du bord de la voilure et le
diagramme d'Enriques est isomorphe au cordage du cerf-volant. 
\medskip

\medskip 
\begin{figure}[h!] 
\vspace*{6mm}
\labellist \small\hair 2pt 
  \pinlabel{$E_0$} at 290 97
  \pinlabel{$E_1$} at 153 146
  \pinlabel{$E_2$} at 365 152
  \pinlabel{$E_3$} at 160 280
  \pinlabel{$E_4$} at 62 115
  \pinlabel{$E_5$} at 28 242
  \pinlabel{$E_6$} at 270 217
  \pinlabel{$E_7$} at 365 239
  \pinlabel{$E_8$} at 3 186
  \pinlabel{$E_9$} at 104 331
  \pinlabel{$E_{10}$} at 195 350
  \pinlabel{$E_{11}$} at 316 335
  \pinlabel{$E_{12}$} at 428 350
  \pinlabel{$E_{13}$} at 478 276
  \pinlabel{$E_{14}$} at 260 370
  \pinlabel{l'axe} at 310 50
\endlabellist 
\centering 
\includegraphics[scale=0.50]{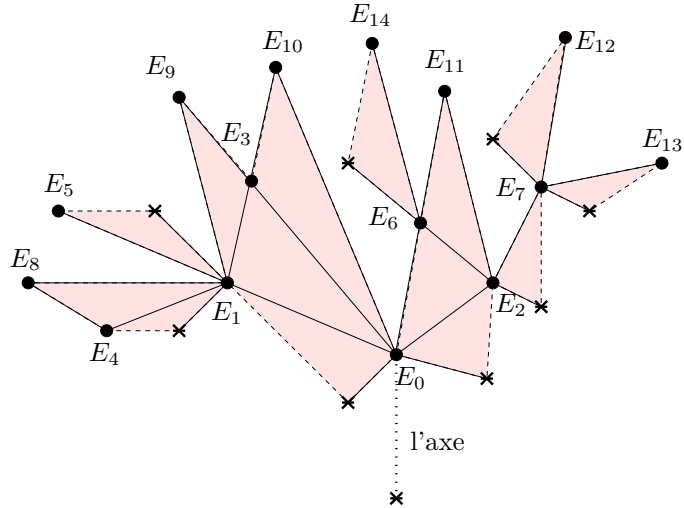} 
\vspace*{0mm} \caption{La voilure de la constellation de la Figure
  \ref{fig:Suitecl}}   
\label{fig:Voilex} 
\end{figure}

Notons par $\mathcal{C} \subset \mathcal{C}_O$ une constellation finie
et par $\mathcal{C}' \supset \mathcal{C}$ une constellation qui 
contient un astre de plus $A \in \mathcal{C}_O$. Expliquons
comment construire la voilure $\mathcal{V}(\mathcal{C}')$ de
$\mathcal{C}'$ {\`a} partir de celle $\mathcal{V}(\mathcal{C})$ de 
$\mathcal{C}$. Le proc{\'e}d{\'e} de construction est exactement le m{\^e}me
que celui d{\'e}crit dans la Section \ref{constvoil} pour les voilures
abstraites. Nous utiliserons la m{\^e}me num{\'e}rotation :

\begin{enumerate}
  \item Si $\mathcal{C}=O$, alors $\mathcal{V}(\{O\})$ est r{\'e}duit
    {\`a} l'axe. On note par $I(O)$ le sommet initial de l'axe et par
    $T(O)$ le sommet terminal.  

  \item Si $\mathcal{C}\neq \emptyset$, alors nous avons deux cas : 

    \begin{enumerate}
        \item \emph{Si $A$ est un astre libre} et que $B= p_D(A)$, on 
          colle une demi-voile $v(A)$ {\`a} $\mathcal{V}(\mathcal{C})$ en
          identifiant son sommet de base {\`a} $T(B)$. On note par
          $I(A)$ le sommet {\'e}toil{\'e} de la demi-voile et par $T(A)$
          son sommet terminal, vus comme points de la voilure obtenue
          apr{\`e}s recollement. 
        
        \item \emph{Si $A$ est un astre satellite} et que $B=p_D(A),
          C=p_I(A)$, on colle une voile simple $v(A)$ {\`a}
          $\mathcal{V}(\mathcal{C})$ en recollant sa base au c{\^o}t{\'e}
          lat{\'e}ral $BC$ de $\mathcal{V}(\mathcal{C})$ par
          l'unique isomorphisme affine qui envoie le
          deuxi{\`e}me sommet de base de la voile simple sur $B$. On note
          par $I(A)$ le point {\'e}toil{\'e} de la nouvelle voile simple
          et par $T(A)$ son sommet terminal, vus comme points de la
          voilure obtenue apr{\`e}s recollement. 
    \end{enumerate}

\end{enumerate}

Comme {\`a} chaque voilure est associ{\'e} canoniquement un cerf-volant, on
obtient aussi \emph{le cerf-volant} $\mathcal{KV}(\mathcal{C})$ \emph{de
  la constellation} $\mathcal{C}$.

\begin{figure}[h!] 
\vspace*{6mm}
\labellist \small\hair 2pt 
  \pinlabel{$E_0$} at 290 97
  \pinlabel{$E_1$} at 153 146
  \pinlabel{$E_2$} at 365 152
  \pinlabel{$E_3$} at 160 280
  \pinlabel{$E_4$} at 62 115
  \pinlabel{$E_5$} at 28 242
  \pinlabel{$E_6$} at 275 217
  \pinlabel{$E_7$} at 370 239
  \pinlabel{$E_8$} at 3 186
  \pinlabel{$E_9$} at 104 331
  \pinlabel{$E_{10}$} at 195 350
  \pinlabel{$E_{11}$} at 316 335
  \pinlabel{$E_{12}$} at 428 350
  \pinlabel{$E_{13}$} at 478 276
  \pinlabel{$E_{14}$} at 261 370

  \pinlabel{${\bf O}$} at 281 -5
  \pinlabel{${\bf A_1}$} at 257 69
  \pinlabel{${\bf A_2}$} at 364 86
  \pinlabel{${\bf A_3}$} at 227 148
  \pinlabel{${\bf A_4}$} at 108 113
  \pinlabel{${\bf A_5}$} at 96 243
  \pinlabel{${\bf A_6}$} at 300 163
  \pinlabel{${\bf A_7}$} at 408 142
  \pinlabel{${\bf A_8}$} at 83 158
  \pinlabel{${\bf A_9}$} at 150 227
  \pinlabel{${\bf A_{10}}$} at 219 208
  \pinlabel{${\bf A_{11}}$} at 325 210
  \pinlabel{${\bf A_{12}}$} at 370 276
  \pinlabel{${\bf A_{13}}$} at 448 219
  \pinlabel{${\bf A_{14}}$} at 228 275
\endlabellist 
\centering 
\includegraphics[scale=0.60]{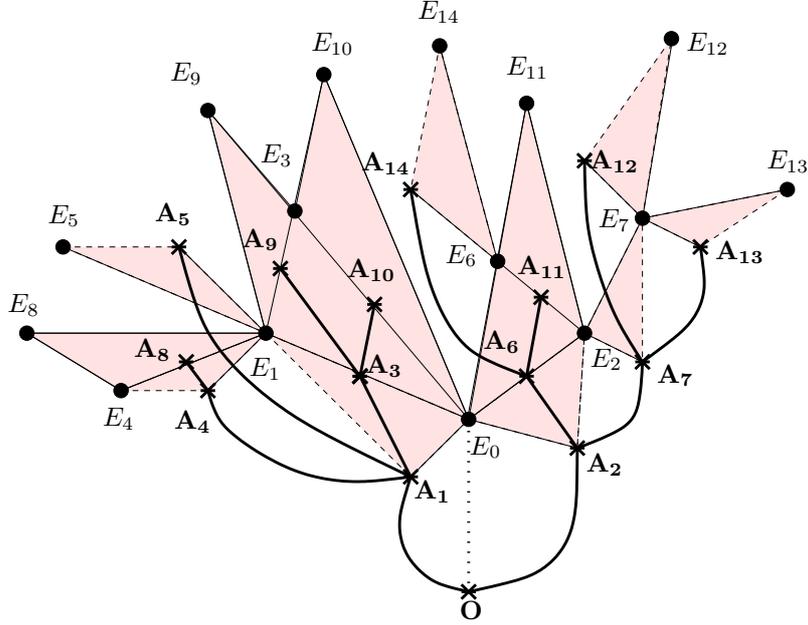} 
\vspace*{0mm} \caption{Le cerf-volant de la constellation de la Figure
  \ref{fig:Suitecl}}   
\label{fig:Cerfconsex} 
\end{figure}

\begin{example}
Dans la Figure \ref{fig:Voilex} est repr{\'e}sent{\'e}e la voilure
de la constellation dont la suite associ{\'e}e
d'{\'e}clatements a {\'e}t{\'e} sch{\'e}matis{\'e}e dans la Figure \ref{fig:Suitecl}. 
Dans la Figure
\ref{fig:Cerfconsex} est repr{\'e}sent{\'e} le cerf-volant associ{\'e}. 
\end{example}

Le th{\'e}or{\`e}me suivant  explique comment retrouver le
diagramme d'Enriques et le graphe dual d'une constellation {\`a} partir
du cerf-volant associ{\'e}. Sa preuve est r{\'e}cursive, en regardant ce qui
se passe {\`a} chaque {\'e}tape d'{\'e}clatement.

\begin{theorem} \label{isos} Soit $\mathcal{C}$ une constellation
  finie centr{\'e}e en $O$.  
   \begin{enumerate}
      \item Le diagramme d'Enriques $\mathcal{E}(\mathcal{C})$ est 
                isomorphe au cordage  
                du cerf-volant $\mathcal{KV}(\mathcal{C})$ de la
                constellation par un isomorphisme qui envoie chaque
                astre $A$ de $\mathcal{C}$ dans $I(A)$. Les c{\^o}t{\'e}s courbes  
                de $\mathcal{E}(\mathcal{C})$ correspondent aux cordes
                libres de 
                $\mathcal{KV}(\mathcal{C})$.  
                Les segments droits maximaux de $\mathcal{E}(\mathcal{C})$
                correspondent aux segments  
                droits maximaux dans la r{\'e}alisation affine canonique
                des voiles compl{\`e}tes de la voilure 
                $\mathcal{V}(\mathcal{C})$. 
      
      \item Le graphe dual $\mathcal{D}(\mathcal{C})$ est isomorphe au graphe
        obtenu comme  
                union des c{\^o}t{\'e}s lat{\'e}raux de la voilure
                $\mathcal{V}(\mathcal{C})$ par un isomorphisme qui
                envoie chaque astre $A$ de $\mathcal{C}$ dans
                $T(A)$. L'auto-intersection d'une 
                composante du diviseur  
                exceptionnel correspondant {\`a} un sommet plein $v$ de
                $\mathcal{KV}(\mathcal{C})$  
                est {\'e}gale {\`a} l'oppos{\'e} du nombre de voiles {\'e}l{\'e}mentaires
                arrivant en $v$, l'axe y compris. 
   \end{enumerate}
\end{theorem}

Remarquons aussi que deux sommets pleins de
$\mathcal{KV}(\mathcal{C})$ sont reli{\'e}s par une ar{\^e}te si et
seulement si, lors du processus d'{\'e}clatement des astres de
$\mathcal{C}$, on trouve un mod{\`e}le sur lequel les centres des deux
valuations divisorielles associ{\'e}es se rencontrent.

\begin{figure}[h!] 
\vspace*{6mm}
\labellist \small\hair 2pt 
  \pinlabel{$E_0$} at 310 97
  \pinlabel{$E_1$} at 170 146
  \pinlabel{$E_2$} at 375 152
  \pinlabel{$E_3$} at 200 249
  \pinlabel{$E_4$} at 62 115
  \pinlabel{$E_5$} at 35 242
  \pinlabel{$E_6$} at 285 207
  \pinlabel{$E_7$} at 375 239
  \pinlabel{$E_8$} at 3 186
  \pinlabel{$E_9$} at 120 331
  \pinlabel{$E_{10}$} at 195 355
  \pinlabel{$E_{11}$} at 325 335
  \pinlabel{$E_{12}$} at 430 370
  \pinlabel{$E_{13}$} at 488 280
   \pinlabel{$E_{14}$} at 270 375
\endlabellist 
\centering 
\includegraphics[scale=0.50]{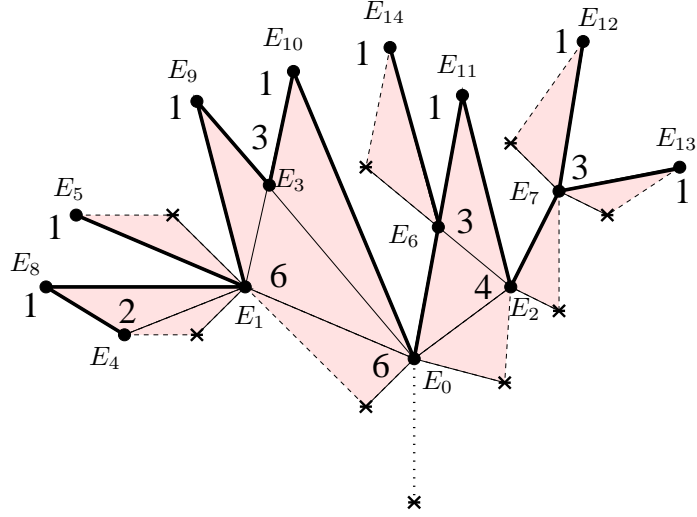} 
\vspace*{0mm} \caption{Plongement du graphe dual dans la voilure  de
  la constellation de la Figure \ref{fig:Suitecl}}   
\label{fig:Dualbord} 
\end{figure}

\begin{example}
Pour l'exemple r{\'e}current de cet article, le diagramme
d'Enriques est visible sur la Figure \ref{fig:Cerfconsex} : c'est le cordage du
cerf-volant. Les arcs courbes et droits sont visibles
directement. Pour d{\'e}terminer les segments qui vont tout droit, on
utilise la Proposition \ref{toutdroit}. On voit alors que la ligne
polygonale $A_1 A_3 A_{10}$ est droite, mais que $A_2 A_6 A_{11}$ ne
l'est pas, ce qui est conforme {\`a} la Figure \ref{fig:Enriqex}. Quant
au graphe dual, nous l'avons repr{\'e}sent{\'e} en traits gras sur la
Figure \ref{fig:Dualbord}. {\`A} c{\^o}t{\'e} de chaque sommet est 
{\'e}crit le nombre de voiles {\'e}l{\'e}mentaires y aboutissant, l'axe y
compris. On v{\'e}rifie ainsi que l'on obtient bien le graphe de la
Figure \ref{fig:Dualex}.
\end{example}

\medskip 

Avoir plong{\'e} {\`a} la fois le diagramme d'Enriques et le graphe dual
dans le m{\^e}me espace de mani{\`e}re {\`a} lire localement dans
cet espace leurs structures suppl{\'e}mentaires (ar{\^e}tes allant tout
droit, auto-intersections) permet d'enrichir la compr{\'e}hension de
tout algorithme de passage de l'un {\`a} l'autre. En effet, {\'e}tant
donn{\'e}e une partie de l'un des graphes, on peut d{\'e}terminer ainsi de
quelle partie de l'autre graphe d{\'e}pend sa structure.

\section{Interpr{\'e}tation valuative} \label{intval}

Les voilures correspondant {\`a} toutes les constellations finies
centr{\'e}es en $O$ peuvent  
{\^e}tre canoniquement recoll{\'e}es. L'on obtient alors la voilure
$\mathcal{V}(\mathcal{C}_O)$  
du firmament $\mathcal{C}_O$ tout entier. On peut 
lui donner une interpr{\'e}tation valuative, analogue {\`a} celle de l'arbre
valuatif de Favre et  
Jonsson \cite{FJ 04}. De plus, il y a une mani{\`e}re naturelle de projectifier 
$\mathcal{V}(\mathcal{C}_O)$ pour obtenir cet arbre valuatif. C'est
ce que nous allons voir dans cette section. 
\medskip

Dans ce qui suit, pour abr{\'e}ger nous noterons $\mathcal{O}:=
\mathcal{O}_{S,O}$, $F:= F_{S,O}$. Soit $\mathcal{M}$ l'id{\'e}al
maximal de $\mathcal{O}$. 

\begin{definition}
  {\'E}tendons la relation d'ordre usuelle de $\R$ {\`a} $\R\cup
  \{\infty\}$ en posant $\infty > \lambda$, pour tout $\lambda \in
  \R$.   Une {\bf valuation de $F$ dominant $O$} est une fonction 
   $\nu:  F \rightarrow \R_+ \cup \{\infty\}$ telle que :
   
 \begin{enumerate}
    \item $\nu(xy)=\nu(x) + \nu(y)$ pour tous $x,y \in F$ ; 

    \item $\nu(x+y) \geq \min(\nu(x), \nu(y))$  pour tous $x,y \in F$ ;

    \item $\nu(\lambda):= \left\{ \begin{array}{ll}
                                    0 & \mbox{si } \lambda \in \C^* \\
                                    \infty & \mbox{si } \lambda=0
                                  \end{array} \right. ;$

    \item $\nu(\mathcal{M}) \subset \R_+^* \cup \{\infty\}$.
 \end{enumerate} 
\end{definition}

Notons par $\mathcal{V}_{S,O}$ l'ensemble des valuations de $F$
dominant $O$ et par $\mathcal{A}_{S,O}$ le sous-ensemble de
$\mathcal{V}_{S,O}$ des valuations normalis{\'e}es par la condition :
 \begin{equation} \label{norm}
   \min \nu(\mathcal{M}) =1.  
 \end{equation} 
Comme l'expliquent Favre et Jonsson de mani{\`e}re d{\'e}taill{\'e}e dans
\cite{FJ 04}, l'ensemble $\mathcal{V}_{S,O}$ admet une topologie
naturelle d'espace fonctionnel localement compact, qui fait du
sous-espace topologique $\mathcal{A}_{S,O}$ un \emph{arbre r{\'e}el}
compact (d'o{\`u} la 
notation $\mathcal{A}$ pour le d{\'e}signer).

Si $A\in \mathcal{C}_O$, notons par $\nu_A$ la  valuation divisorielle
associ{\'e}e. Elle peut {\^e}tre d{\'e}finie des deux mani{\`e}res
{\'e}quivalentes suivantes, en partant d'un mod{\`e}le $(\Sigma,E)
\stackrel{\pi}{\rightarrow} (S,O)$ contenant $A$ :

$\bullet$ si $f \in F$, alors $\nu_A(f)$ est la multiplicit{\'e} de la
fonction $f\circ \pi$ au point $A$ ; 

$\bullet$ si $f \in F$, alors $\nu_A(f)$ est l'ordre d'annulation de
$f \circ \pi \circ \pi_A$ le long de $E_A$, o{\`u} $\Sigma_A
\stackrel{\pi_A}{\rightarrow} \Sigma$ est l'{\'e}clatement de $A$ dans $\Sigma$ et
que $E_A$ est le diviseur exceptionnel ainsi cr{\'e}{\'e}. 

\medskip

Reprenons les notations de la section
pr{\'e}c{\'e}dente : $\mathcal{C}$ d{\'e}signe donc une constellation finie
et $A$ est un astre que l'on rajoute {\`a} $\mathcal{C}$. Nous
expliquons {\`a} pr{\'e}sent comment plonger canoniquement la voile
{\'e}l{\'e}mentaire $v(A)$ associ{\'e}e {\`a} $A$ dans l'espace valuatif
$\mathcal{V}_{S,O}$. Nous allons changer l'ordre consid{\'e}r{\'e}
auparavant, en traitant d'abord le cas o{\`u} $A$ est \emph{satellite},
ensuite celui o{\`u} $A$ est \emph{libre} et enfin celui o{\`u} $A=O$. 

\medskip 
\begin{figure}[h!] 
\vspace*{6mm}
\labellist \small\hair 2pt 
  \pinlabel{$E_B=\{x_B=0\}$} at 83 -20
  \pinlabel{$E_C=\{x_C=0\}$} at -50 80
  \pinlabel{$A$} at 96 100
  \pinlabel{$O$} at 483 63
  \pinlabel{$S$} at 570 120
  \pinlabel{$S_{\mathcal{C}}$} at 195 120
  \pinlabel{$\pi_{\mathcal{C}}$} at 296 100
\endlabellist 
\centering 
\includegraphics[scale=0.50]{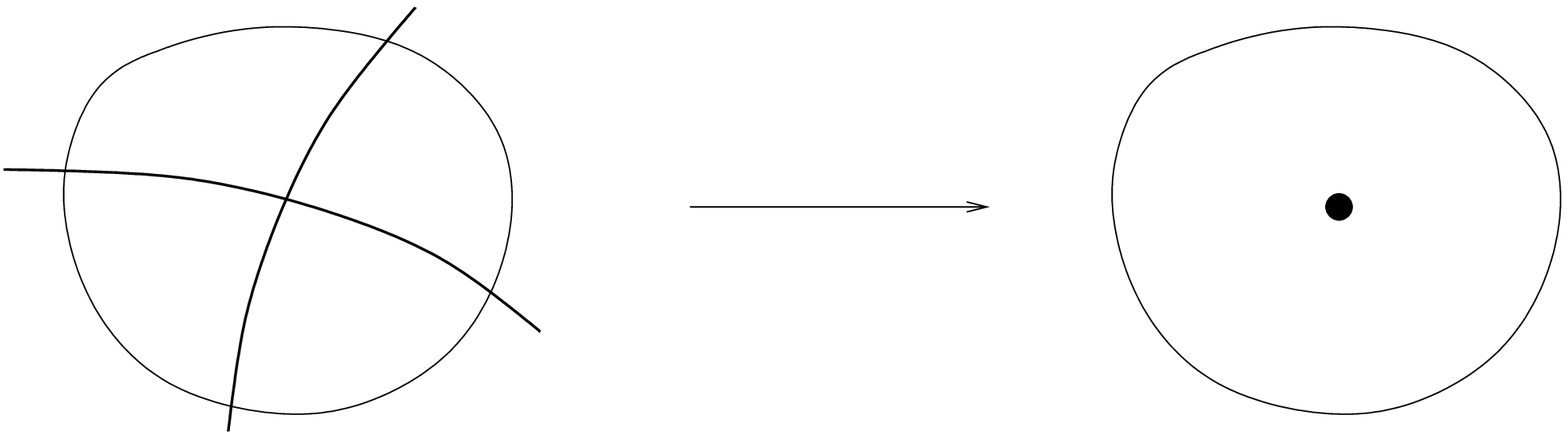} 
\vspace*{5mm} \caption{Le cas d'un astre satellite}   
\label{fig:Sat} 
\end{figure}

$\bullet$ \emph{Supposons que $A$ est satellite} (voir la Figure
\ref{fig:Sat}).  
Avec les notations de l'{\'e}quation (\ref{eclconst}), notons par $E_B$
et $E_C$ les transform{\'e}es strictes sur la surface $S_{\mathcal{C}}$
des diviseurs exceptionnels obtenus en {\'e}clatant les points
$B=p_D(A), C=p_I(A)$. Ces deux courbes s'intersectent transversalement au
point $A$. Notons par $N_A$ le r{\'e}seau abstrait engendr{\'e} par les valuations
divisorielles $\nu_B$ et $\nu_C$. Notons par $(e_B, e_C)$ la base de $N_A$
correspondant au couple $(\nu_B, \nu_C)$. 

\begin{definition}
   Une valuation $\nu\in \mathcal{V}_{S,O}$ est dite {\bf monomiale}
   par rapport {\`a} $\nu_B$ et $\nu_C$ s'il existe $(b,c) \in \R_+^2
   \setminus 0$ 
   tels que pour tout $f \in F^*$ : 
    $$ \nu (f) = \min \{ b\cdot m_B + c \cdot m_C \: | \: c_{m_B, m_C}
    \neq 0, \: f\circ \pi_{\mathcal{C}} = \sum_{(m_B, m_C)}  
       c_{m_B,m_C} x_B^{m_B} x_C^{m_C} \}$$
   o{\`u} $(x_B, x_C)$ est un syst{\`e}me de coordonn{\'e}es locales au
   voisinage de $A\in S_{\mathcal{C}}$ tel que $E_B, E_C$ soient
     d{\'e}finis par les {\'e}quations $x_B=0$, respectivement $x_C=0$. 
  Nous noterons par $b \: \nu_B \oplus c \: \nu_C$ la valuation monomiale
  pr{\'e}c{\'e}dente. 
\end{definition}

Cette d{\'e}finition est ind{\'e}pendante du choix du syst{\`e}me de
coordonn{\'e}es locales. Le nom est motiv{\'e} par le fait que la
valuation d'une fonction se d{\'e}termine uniquement {\`a} partir des
valuations $\nu_B(x_B^{m_B} x_C^{m_C})$ et $\nu_C(x_B^{m_B}
x_C^{m_C})$ des mon{\^o}mes intervenant dans l'{\'e}criture de
$f\circ\pi_{\mathcal{C}}$ dans le syst{\`e}me de coordonn{\'e}es $(x_B, x_C)$.

De cette mani{\`e}re, les notations {\'e}tant celles du d{\'e}but de la
Section \ref{modaff}, le c{\^o}ne $\sigma(e_B, e_C)$ de l'espace vectoriel
$(N_A)_{\R}$  se plonge dans $\mathcal{V}_{S,O}$, en associant {\`a}
chaque vecteur $b\cdot e_B + c\cdot e_C$ la valuation monomiale $b\:
\nu_B \oplus c\: \nu_C$. Notons par $\sigma(\nu_B, \nu_C)$ son
image. On prend comme voile {\'e}l{\'e}mentaire $v(A)$  
le triangle affine du plan $(N_A)_{\R}$ de sommets $e_B, e_C, e_B +
e_C$. Par le plongement pr{\'e}c{\'e}dent, il se r{\'e}alise comme triangle
dans  $\mathcal{V}_{S,O}$ de sommets $\nu_B, \nu_C, \nu_B \oplus
\nu_C$. Le fait que les sous-triangles de l'espace
topologique 
$\mathcal{V}_{S,O}$ correspondant {\`a} $v(A), v(B), v(C)$ se retrouvent
recoll{\'e}s comme d{\'e}crit dans la 
construction de la voilure $\mathcal{V}(\mathcal{C})$ provient du
lemme {\'e}l{\'e}mentaire suivant :

\begin{lemma} \label{sumval}
  On a l'{\'e}galit{\'e} suivante de valuations :
   $\nu_A= \nu_B \oplus \nu_C.$
\end{lemma}

\medskip
$\bullet$ \emph{Supposons que $A$ est libre}. Notons par $E_B$ la
transform{\'e}e stricte sur la surface $S_{\mathcal{C}}$ du diviseur
  exceptionnel obtenu en {\'e}clatant $B= p_D(A)$. C'est l'unique
  composante de $E_{\mathcal{C}}$ qui contient le point $A$. Dans ce
  cas on choisit une \emph{curvette} passant par $A$, c'est-{\`a}-dire
  un germe de courbe lisse transverse {\`a} $E_B$. Cette curvette jouera
  le r{\^o}le de $E_C$. Notons par $\nu_{curv}\in
  \mathcal{V}(\mathcal{C})$ la valuation  divisorielle associ{\'e}e. 

  On fait la m{\^e}me construction que pr{\'e}c{\'e}demment, le couple
  $(\nu_B, \nu_{curv})$ de valuations jouant le m{\^e}me r{\^o}le que
  $(\nu_B, \nu_C)$. Cette fois-ci une partie des valuations monomiales
  $b \: \nu_B \oplus c \: \nu_{curv}$ d{\'e}pendent bien s{\^u}r du choix
  de la curvette, mais le point important est que la moiti{\'e} du
  c{\^o}ne de ces valuations n'en d{\'e}pend pas :

\begin{lemma}
  L'intersection dans $\mathcal{V}_{S,O}$ des c{\^o}nes $\sigma(\nu_B,
  \nu_{curv})$, lorsque la curvette varie, est {\'e}gale au c{\^o}ne
  $\sigma(\nu_B, \nu_A)$.  Plus pr{\'e}cis{\'e}ment, parmi les valuations de
  la forme $b \:
  \nu_B \oplus c \: \nu_{curv}$, celles v{\'e}rifiant $b \geq c$ sont
  exactement les valuations ind{\'e}pendantes 
  du choix de la curvette, et co{\"\i}ncident avec les valuations du
  c{\^o}ne $\sigma(\nu_B, \nu_A)$.
\end{lemma}

Ce lemme d{\'e}coule directement du Lemme \ref{sumval}.

En faisant la construction analogue {\`a} celle du cas pr{\'e}c{\'e}dent, on
voit gr{\^a}ce {\`a} ce lemme que la demi-voile associ{\'e}e {\`a}
l'astre libre $A$ est le  
triangle affine de $\mathcal{V}_{S,O}$ de sommets $\frac{1}{2}\nu_A,
\nu_A, \nu_B$. Le sommet {\'e}toil{\'e} correspond {\`a} $\frac{1}{2}\nu_A$,
le sommet terminal {\`a} $\nu_A$ et le sommet de base {\`a} $\nu_B$.
\medskip

$\bullet$ \emph{Supposons que $A=O$}. Consid{\'e}rons un syst{\`e}me
quelconque de coordonn{\'e}es centr{\'e} en $O$. Il lui correspond un
triangle affine, p{\'e}tale initial du lotus associ{\'e}. Lorsque l'on
varie le syst{\`e}me de coordonn{\'e}es, l'intersection
de tous ces p{\'e}tales est {\'e}gale au segment affine des valuations de
la forme $\lambda \nu_O$, avec $\lambda \in [\frac{1}{2}, 1]$. On
envoie bijectivement l'axe sur ce segment par l'unique application
affine qui associe $\frac{1}{2} \nu_A$ au sommet initial  $O$ et
$\nu_A$ au sommet terminal $E_O$.

\medskip

En conclusion des consid{\'e}rations pr{\'e}c{\'e}dentes :

\begin{proposition}
  La voilure d'une constellation se plonge canoniquement dans
  $\mathcal{V}_{S,O}$ en pr{\'e}servant les structures affines des voiles
  compl\`etes.  
\end{proposition}

On peut montrer que le cerf-volant se plonge aussi canoniquement dans
le m{\^e}me espace. L'id{\'e}e est de consid{\'e}rer pour chaque 
astre libre $A$ contenue sur une composante $E_i$ une curvette passant par
$A$ et un syst{\`e}me de coordonn{\'e}es centr{\'e} en l'astre $A_i$
d'ont l'un des axes de coordonn{\'e}es soit l'image de cette
curvette. On relie les valuations repr{\'e}sentatives de $A$ et $A_i$
par un segment dans le plan des valuations monomiales par rapport {\`a}
ce syst{\`e}me de coordonn{\'e}es, et on montre que ce segment est
ind{\'e}pendant des choix. 

Dans la suite de la section, pour chaque astre $A$, nous identifierons
$v(A)$ {\`a} un 
sous-triangle de $\mathcal{V}_{S,O}$. Il est plong{\'e} dans l'espace
des valuations monomiales par rapport {\`a} $(\nu_B, \nu_C)$ (si $A$ est
satellite) ou $(\nu_A, \nu_B)$ (si $A$ est libre). Tra{\c c}ons sur ce
triangle le feuilletage $\mathcal{F}_A$ obtenu en l'intersectant avec
les droites partant de l'origine dans le c{\^o}ne des valuations monomiales
correspondant. L'espace des feuilles s'identifie canoniquement (par
intersection) {\`a} l'union des c{\^o}t{\'e}s lat{\'e}raux de $v(A)$. Mais ces
c{\^o}t{\'e}s lat{\'e}raux sont des valuations normalis{\'e}es par la condition
(\ref{norm}). En recollant 
les feuilletages des voiles {\'e}l{\'e}mentaires de la voilure
$\mathcal{V}(\mathcal{C})$, on obtient un feuilletage
$\mathcal{F}(\mathcal{C})$. Les applications de passage au
quotient par les feuilles se recollent en une application :
$$
 \displaystyle{ \mathcal{V}(\mathcal{C})
   \stackrel{\phi_{\mathcal{C}}}{\longrightarrow} 
 \mathcal{D}(\mathcal{C})  }
$$
dans laquelle la voilure $\mathcal{V}(\mathcal{C})$ est vue comme
sous-espace de l'espace des valuations $\mathcal{V}_{S,O}$ et le
graphe dual $\mathcal{D}(\mathcal{C})$ est vu comme sous-espace de
l'arbre $\mathcal{A}(\mathcal{C})$ des valuations normalis{\'e}es. 

Si on a une inclusion $\mathcal{C} \subset \mathcal{C}'$ de
constellations finies, on a des r{\'e}tractions naturelles
$\mathcal{V}(\mathcal{C}') \rightarrow \mathcal{V}(\mathcal{C})$ et  
$\mathcal{D}(\mathcal{C}') \rightarrow \mathcal{D}(\mathcal{C})$ 
telles que le diagramme suivant soit commutatif :

$$\xymatrix{
       \mathcal{V}(\mathcal{C}') \ar[d] \ar[r]^{\phi_{\mathcal{C}'}}
      &  \mathcal{D}(\mathcal{C}')   \ar[d]  
       \\
        \mathcal{V}(\mathcal{C})  \ar[r]^{\phi_{\mathcal{C}}}
      &  \mathcal{D}(\mathcal{C})   }$$
On peut prendre alors les limites projectives des morphismes
$\phi_{\mathcal{C}}$. Notons par :
  $$\mathcal{V}(\mathcal{C}_O):=
  \lim_{\longleftarrow}\mathcal{V}(\mathcal{C})  $$ 
\emph{la voilure du firmament}. Comme la limite projective des graphes
duaux $\mathcal{D}(\mathcal{C})$ s'identifie {\`a} l'arbre valuatif
$\mathcal{A}(\mathcal{C})$ (voir \cite{FJ 04}), on obtient :

\begin{proposition}
 L'application de passage au quotient de la voilure
 $\mathcal{V}(\mathcal{C}_O)$ du firmament par
 le feuilletage limite projective des feuilletages
 $\mathcal{F}(\mathcal{C})$ s'identifie naturellement {\`a}
 l'arbre valuatif $\mathcal{A}_{S,O}$. 
\end{proposition}

\section{Le lotus et les fractions continues}
\label{seclot}

Dans cette section j'indique comment
interpr{\'e}ter g{\'e}om{\'e}triquement les d{\'e}veloppements en fractions
continues des nombres positifs {\`a} l'aide du lotus. 
 Cette interpr{\'e}tation est indispensable 
d{\`e}s qu'on veut d{\'e}crire le cerf-volant d'un germe de courbe plane en
termes des exposants de Newton-Puiseux de ses branches (voir \cite{GP
  09}). 
\medskip

Comme dans \cite{PP 07}, j'utiliserai les notations suivantes pour les
fractions continues \emph{euclidiennes} (uniquement des signes $+$)  et
\emph{de Hirzebruch-Jung} (uniquement des signes $-$) : 
 $$ [x_1,x_2, ...]^{\pm} :=  x_1 \pm \cfrac{1}{x_2 \pm
                          \cfrac{1}{\cdots}}\: . $$

Consid{\'e}rons {\`a} nouveau un r{\'e}seau bidimensionnel $N$, muni d'une
base $(e_1, e_2)$. 
Notons par $D(e_1, e_2)$ la droite de $N_{\R}$ qui joint les points $e_1,
e_2$. Les couples de points entiers successifs sur cette droite sont
de la forme $((1-a)e_1 + a e_2, -a e_1 + (a+1) e_2)$, o{\`u} $a \in
\Z$. Ils forment des bases de $N$. 

Consid{\'e}rons {\`a} pr{\'e}sent la droite passant par $0$ et parall{\`e}le
{\`a} la droite $D(e_1, e_2)$ : 
   $$D_{\infty}(e_1, e_2):= \R (e_1
     -e_2) \subset N_{\R}.$$ 
Notons par $P(e_1, e_2)$ le demi-plan ouvert bord{\'e} par
  $D_{\infty}(e_1, e_2)$ et contenant $D(e_1, e_2)$, puis par
  $\Delta_0(e_1, e_2)$ le triangle de sommets $0, e_1, e_2$. On a
  {\'e}videmment : 

\begin{lemma}
  L'union des c{\^o}nes $\sigma((1-a)e_1 + a e_2, -a e_1 + (a+1) e_2)$
  pour $a \in \Z$ est {\'e}gale {\`a} $P(e_1, e_2)\cup 0$.  
\end{lemma} 

\begin{figure}[h!] 
\vspace*{5mm}
\labellist \small\hair 2pt 
  \pinlabel{$0$} at 130 130
  \pinlabel{$e_1$} at  274 140
  \pinlabel{$e_2$} at  150 270
  \pinlabel{$e_1+ e_2$} at  290 277
  \pinlabel{$2e_1 - e_2$} at  390 10
  \pinlabel{$-e_1 + 2 e_2$} at  10 375
  \pinlabel{$2 e_1 + e_2$} at  370 250
  \pinlabel{$ e_1 + 2 e_2$} at  256 360
  \pinlabel{$D_{\infty}(e_1, e_2)$} at 235 15
  \pinlabel{$D(e_1, e_2)$} at 5 320
\endlabellist 
\centering 
\includegraphics[scale=0.7]{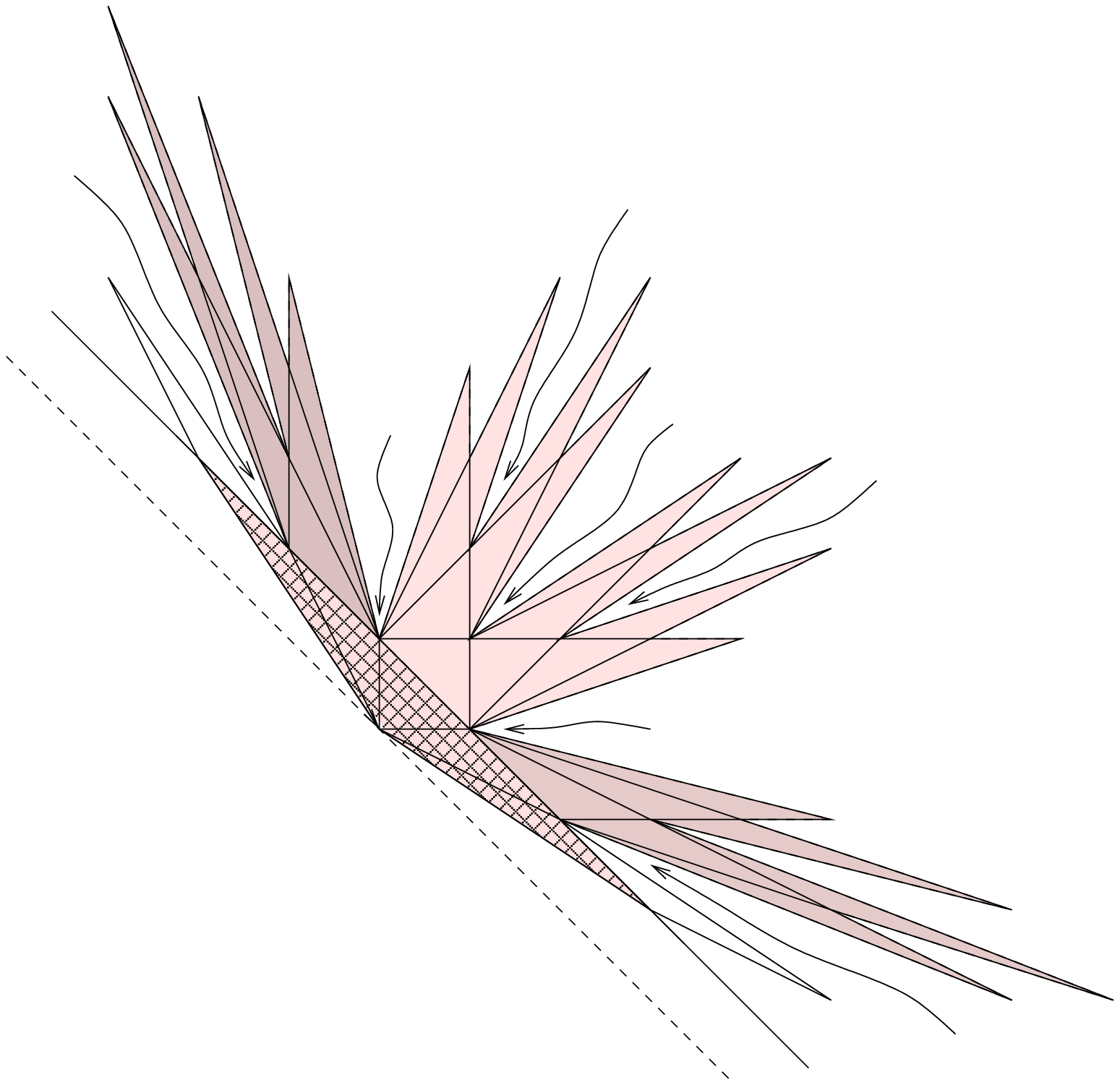} 
\vspace*{5mm} \caption{Le grand lotus $\overline{\mathcal{L}}(e_1, e_2)$} 
\label{fig:Grandlotus} 
\end{figure} 

Introduisons alors :

\begin{definition} {\bf Le grand lotus} $\overline{\mathcal{L}}(e_1,
e_2)$ associ{\'e} {\`a} la base $(e_1, e_2)$ est l'union des lotus
associ{\'e}s {\`a} tous les c{\^o}nes du lemme pr{\'e}c{\'e}dent, ainsi que des triangles
$\Delta_0((1-a)e_1 + a e_2, -a e_1 + (a+1) e_2)$ (voir la Figure
\ref{fig:Grandlotus}) :
  $$ \overline{\mathcal{L}}(e_1, e_2) := 
     \bigcup_{a \in \Z} [\mathcal{L}((1-a)e_1 + a e_2, -a e_1 + (a+1)
     e_2)\cup \Delta_0((1-a)e_1 + a e_2, -a e_1 + (a+1) e_2) ].$$
\end{definition}

Les sommets du grand lotus forment l'ensemble $\mathrm{Prim} (P(e_1,
e_2))$ des vecteurs primitifs du demi-plan ouvert $P(e_1, e_2)$,
auquel on rajoute $0$.

Consid{\'e}rons {\`a} pr{\'e}sent $\mathbb{H}(N_{\R})$, le plan hyperbolique  
dont l'horizon est la droite projective r{\'e}elle
$\mathbb{P}(N_{\R})$. 

\begin{remark} Une construction canonique de
$\overline{\mathbb{H}}(N_{\R}):= \mathbb{H}(N_{\R}) \cup
\mathbb{P}(N_{\R})$ peut {\^e}tre 
faite de la mani{\`e}re suivante. A chaque structure presque complexe
$J$ sur $N_{\R}$ (c'est-{\`a}-dire un endomorphisme de $N_{\R}$
v{\'e}rifiant $J^2 =-I$) on associe la d{\'e}composition $N_{\C}=N_J^i
\oplus N_J^{-i}$ en somme directe des espaces propres du complexifi{\'e}
$J_{\C}:N_{\C} \rightarrow N_{\C}$ de $J$. Les droites complexes
$N_J^i$ et $N_J^{-i}$ sont conjugu{\'e}es par rapport {\`a} la conjugaison
canonique $u + iv \rightarrow u -i v$ de $N_{\C}= N_{\R} +
iN_{\R}$. L'application $J \rightarrow N_J^i$ identifie
bijectivement l'ensemble des structures presque complexes sur $N_{\R}$
avec l'ensemble des droites complexes imaginaires de $N_{C}$. Ce
dernier ensemble s'identifie {\`a} 
  $\mathbb{P}(N_{C}) \setminus \mathbb{P}(N_{\R})$, 
c'est-{\`a}-dire au compl{\'e}mentaire d'un cercle dans la sph{\`e}re de
Riemann     $\mathbb{P}(N_{C})$. Chacun des deux h{\'e}misph{\`e}res ainsi
d{\'e}limit{\'e}s repr{\'e}sente les structures presque complexes
d{\'e}finissant l'une des deux orientations de $N_{\R}$. On consid{\`e}re
alors sur chacun des h{\'e}misph{\`e}res l'unique m{\'e}trique hyperbolique
d{\'e}finissant la m{\^e}me structure conforme. La conjugaison restreinte
{\`a} $\mathbb{P}(N_{C}) \setminus \mathbb{P}(N_{R})$ est une isom{\'e}trie renversant
l'orientation h{\'e}rit{\'e}e de celle de $\mathbb{P}(N_{C})$. On peut
donc d{\'e}finir $\mathbb{H}(N_{\R})$ au choix, soit comme 
\emph{l'espace des $J$ pr{\'e}servant une orientation fix{\'e}e de} $N_{\R}$, soit
comme \emph{l'espace des couples non-ordonn{\'e}s} $\{J, -J\}$.
\end{remark}

\begin{figure}[h!] 
\vspace*{5mm}
\labellist \small\hair 2pt 
  \pinlabel{$[e_1 -e_2]$} at 255 -20
  \pinlabel{$[2e_1-e_2]$} at 355 0
  \pinlabel{$[e_1]$} at 487 115
  \pinlabel{$[2e_1 +e_2]$} at 510 395

  \pinlabel{$[e_1 + e_2]$} at 252 520
  \pinlabel{$[e_1+ 2 e_2]$} at -10 385
  \pinlabel{$[e_2]$} at 14 115
  \pinlabel{$[-e_1 + 2 e_2]$} at 140 0
\endlabellist 
\centering 
\includegraphics[scale=0.50]{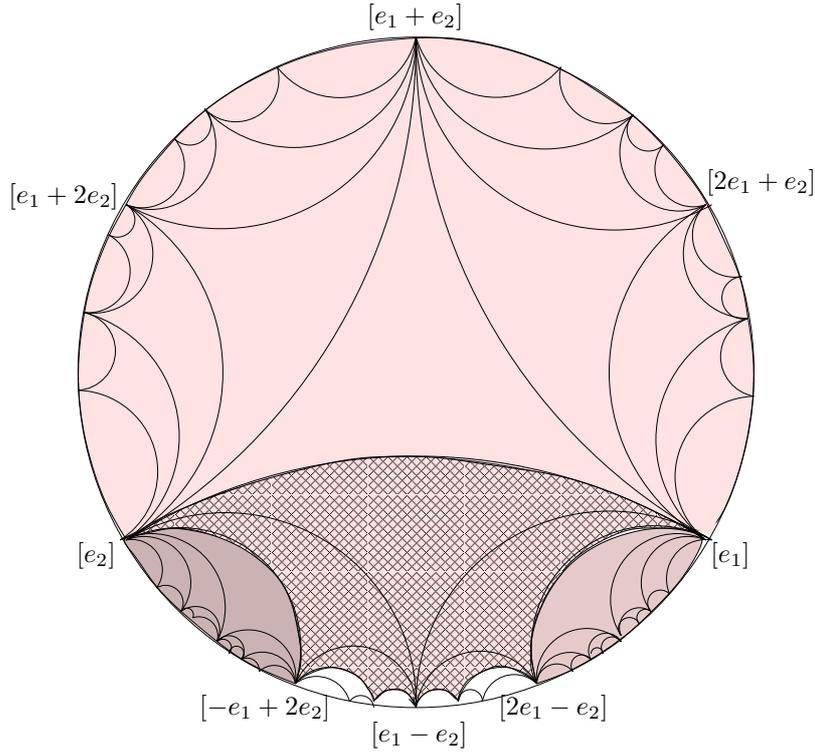} 
\vspace*{5mm} \caption{La triangulation modulaire $\mathcal{M}(N)$} 
\label{fig:Hyperb} 
\end{figure}

Pour chaque $v \in N_{\R} \setminus 
0$, nous noterons par $[v]\in \mathbb{P}(N_{\R})$ le point
correspondant {\`a} la droite $\R v$.  Pour chaque base $(u,v)$ de $N$,
soit $L([u],[v])$ l'unique droite hyperbolique de $H(N_{\R})$ qui
joint $[u]$ et 
$[v]$. Lorsque $(u,v)$ varie parmi toutes les bases de $N$, on
obtient des droites qui forment les ar{\^e}tes d'une triangulation de
$\overline{H}(N_{\R}):= H(N_{\R}) \cup \mathbb{P}(N_{\R})$, appel{\'e}e
\emph{la triangulation modulaire} $\mathcal{M}(N)$. Le grand lotus est
un plongement 
combinatoire de la triangulation modulaire dans $N_{\R}$, associ{\'e}e
canoniquement {\`a} la base $(e_1, e_2)$ :

\begin{proposition} \label{reamod}
   L'application :
   $$ \begin{array}{lccc}
      \Phi : & \mathrm{Prim}( P(e_1, e_2)) \cup 0 & \rightarrow &
      \mathbb{P}(N)\\ 
           &  v \neq o& \rightarrow & [v]\\
           &  0 & \rightarrow & [e_1 -e_2]
              \end{array}  $$
 est une bijection qui envoie le grand
 lotus $\overline{\mathcal{L}}(e_1, e_2)$  en la triangulation modulaire
 $\mathcal{M}(N)$.  
\end{proposition}

Dans la Figure \ref{fig:Hyperb} est repr{\'e}sent{\'e}e la triangulation
modulaire. Dans les Figures \ref{fig:Grandlotus} et \ref{fig:Hyperb}
sont repr{\'e}sent{\'e}es avec les m{\^e}mes couleurs certains
sous-complexes de  $\overline{\mathcal{L}}(e_1, e_2)$ et leurs images par
$\Phi$. En particulier, on voit que le lotus $\mathcal{L}(e_1, e_2)$
correspond {\`a} la partie de la triangulation modulaire situ{\'e}e dans
le demi-plan hyperbolique bord{\'e} par la droite $L([e_1], [e_2])$ et
contenant $[e_1 + e_2]$ dans son adh{\'e}rence.

\begin{remark}
   Dans certains travaux (voir entre autres \cite{S 82} ,\cite{S 85},
   \cite{HT 85}, \cite{H 00}), les
   fractions continues sont interpr{\'e}t{\'e}es g{\'e}om{\'e}triquement {\`a}
   l'aide de la triangulation modulaire $\mathcal{M}(N)$. Par la
   proposition pr{\'e}c{\'e}dente, ces interpr{\'e}tations peuvent se faire
   de mani{\`e}re {\'e}quivalente sur le grand lotus. Dans la suite de cette
   section, 
   j'explique une interpr{\'e}tation diff{\'e}rente, faite uniquement {\`a} l'aide du
   lotus, et ayant l'avantage de se g{\'e}n{\'e}raliser en dimensions plus
   grandes (voir la Section \ref{lotarb}). C'est une
   r{\'e}interpr{\'e}tation de celle  
   de Klein (voir \cite{K 96} et \cite{PP 07}), mais la
   g{\'e}n{\'e}ralisation qu'elle 
   sugg{\`e}re en dimension plus grande est diff{\'e}rente de celle
   propos{\'e}e par Klein et reprise entre autres par Arnold \cite{A 98}. 
\end{remark}

Consid{\'e}rons une demi-droite $l\subset N_{\R}$ d'origine $0$ et
contenue dans l'int{\'e}rieur du c{\^o}ne $\sigma(e_1, e_2)$. Notons par
$G(l)$ l'union des p{\'e}tales du lotus dont l'int{\'e}rieur intersecte
$l$. Nous dirons que $G(l)$ est \emph{la gaine} de $l$. On peut
l'imaginer construite successivement en rajoutant des p{\'e}tales
$\tau_1, \tau_2,...$ {\`a} $\tau_0:= \tau(e_1, e_2)$, au fur et {\`a}
mesure que l'on s'{\'e}loigne de $0$ le long de $l$ : chaque fois que
l'on entre dans un nouveau p{\'e}tale, on le rajoute {\`a} la
suite d{\'e}j{\`a} construite. Consid{\'e}rons deux cas, suivant que $l$ est
rationnelle ou non.

$ \bullet$ 
Si $l$ est \emph{rationnelle}, c'est-{\`a}-dire que $l$ contient des
{\'e}l{\'e}ments non-nuls du r{\'e}seau $N$, notons par $S(l)$ l'unique
{\'e}l{\'e}ment primitif de $N$ contenu dans $l$. Le segment $[0, S(l)]$
ne contient pas d'autre points de $N$ {\`a} part ses extr{\'e}mit{\'e}s, ce
qui permet de montrer que :
 $$[0, S(l)] = G(l) \cap l.$$
Dans ce cas, $G(l)$ contient un nombre fini de p{\'e}tales
$\tau_0,...,\tau_n$ et $S(l)$ est un sommet de $\tau_n$. Pour $i=1,2$,
notons par
$P_i(l)$ la ligne polygonale joignant $e_i$ {\`a} $S(l)$ et contenue
dans le bord de $G(l)$.

\begin{figure}[h!] 
\vspace*{5mm}
\labellist \small\hair 2pt 
  \pinlabel{0} at -5 -5
  \pinlabel{$l$} at 542 408
  \pinlabel{$S(l)$} at 535 340
  \pinlabel{$e_1$} at 75 -10
  \pinlabel{$e_2$} at -10 75
  \pinlabel{$d_1$} at 325 -15
  \pinlabel{$d_2$} at -15 250
  \pinlabel{$P_1(l)$} at 215 110
  \pinlabel{$P_2(l)$} at 100 130
\endlabellist 
\centering 
\includegraphics[scale=0.40]{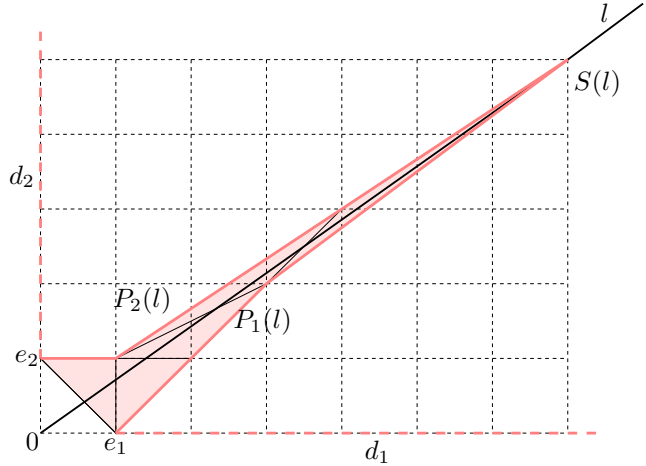} 
\vspace*{5mm} \caption{La gaine de la demi-droite rationnelle
  $\R_+(7e_1 + 5e_2)$} 
\label{fig:Gaine} 
\end{figure}

$\bullet$ 
Si $l$ est \emph{irrationnelle}, la suite $\tau_0, \tau_1,...$
est infinie. Pour $i=1,2$, notons par $P_1(l)$ la ligne polygonale
contenue dans le bord de $G(l)$, partant de $e_i$ et allant vers
l'infini en s'{\'e}loignant toujours strictement de $0$ (ou, de
mani{\`e}re {\'e}quivalente, la composante connexe de $\partial G(l)
\setminus ]e_1, e_2[$ contenant $e_i$). 

Notons par $d_i$ la demi-droite ferm{\'e}e contenue dans la droite
engendr{\'e}e par le vecteur $e_i$, d'origine le point $e_i$ et ne
contenant pas $0$. Posons aussi :
$$ Q_i(l) := \left\{ 
      \begin{array}{ll}
          d_i \cup P_i(l)\cup (l \setminus [0, S(l)]), &
                \mbox{ si $l$ est rationnelle   } \\
           d_i \cup P_i(l), &
                \mbox{ si $l$ est irrationnelle   } 
      \end{array}   \right. $$

\medskip 
Associons aussi {\`a} chaque p{\'e}tale $\tau_i$ l'un des symboles
`$\delta$' ou `$\gamma$', suivant que $\tau_{i+1}$ est attach{\'e} du
c{\^o}t{\'e} \emph{droit} ou \emph{gauche} de $\tau_i$ et convenons que le
symbole du dernier p{\'e}tale (si $l$ est rationnel) est le m{\^e}me que
celui de l'avant-dernier.

\begin{example} \label{exgaine}
  Dans la Figure \ref{fig:Gaine} est repr{\'e}sent{\'e}e la gaine de la
  demi-droite $l$ passant par $S(l)= 7e_1 + 5 e_2$. Elle est form{\'e}e
  de 5 p{\'e}tales. $P_1(l)$ joint dans l'ordre les points $e_1, 3e_1 + 2e_2, S(l)$
  et $P_2(l)$ joint dans l'ordre les points $e_2,e_1 + e_2, S(l)$. Les 
  symboles des p{\'e}tales sont, dans l'ordre, $\delta, \gamma, \gamma,
  \delta, \delta$. 
\end{example}

La proposition suivante relie la construction de la gaine {\`a} la
construction g{\'e}om{\'e}trique de Klein lui permettant de donner une
interpr{\'e}tation des fractions continues (voir \cite{PP 07}). On en
d{\'e}duit l'interpr{\'e}tation g{\'e}om{\'e}trique annonc{\'e}e des fractions continues
{\`a} l'aide du 
lotus. Rappelons que \emph{l'inclinaison} de $l$ dans la base $(e_1,
e_2)$ est le quotient $x_1 /x_2$, si $x_1 e_1 + x_2 e_2$ est un
vecteur directeur de $l$. 

\begin{proposition} \label{gainefrac}
  \begin{enumerate}
  
   \item Pour chaque $i \in \{1,2\}$, la ligne polygonale $Q_i(l)$ est
     le bord de l'enveloppe convexe de 
  l'ensemble des points du r{\'e}seau $N$ contenus dans le c{\^o}ne
  strictement convexe de c{\^o}t{\'e}s $\R_+e_i$ et $l$.

   \item  Regroupons les p{\'e}tales $\tau_0,\tau_1,...$ de la gaine
     $G(l)$ en paquets
     maximaux de triangles successifs ayant le m{\^e}me symbole. On
     consid{\`e}re que le premier paquet est toujours form{\'e} 
     de p{\'e}tales tournant {\`a} droite. D{\'e}signons par $a_1, a_2,...$
     les cardinaux des paquets successifs. Alors, si $\lambda$
     d{\'e}signe l'inclinaison de $l$ dans la base $(e_1, e_2)$, on a :
     $$ \lambda =  [a_1, a_2, ... ]^+. $$

  \end{enumerate}
\end{proposition}

\begin{remark}
  Dans \cite{PP 07} j'ai expliqu{\'e} que l'on pouvait comprendre la
  dualit{\'e} des enveloppes convexes des points entiers situ{\'e}s dans
  deux c{\^o}nes suppl{\'e}mentaires {\`a} l'aide d'un \emph{diagramme en
    zig-zag}. Dans le cas o{\`u} ces c{\^o}nes sont celui de c{\^o}t{\'e}s
  $\R_+ (e_1 -e_2), l$ et celui de c{\^o}t{\'e}s $l, \R_+ (e_2 -e_1)$, la
  ligne en zig-zag associ{\'e}e obtenue en partant du point $e_1 -
  e_2$ est l'union du segment $[e_1 -e_2, e_2]$ et des segments qui
  s{\'e}parent {\`a} l'int{\'e}rieur de la gaine $G(l)$ les unions de
  p{\'e}tales tournant dans le m{\^e}me sens. On pourra comparer ceci aux
  consid{\'e}rations de \cite{HT 85}.  
\end{remark}

Au d{\'e}but de la section, le grand lotus a {\'e}t{\'e} associ{\'e} {\`a} une
base de $N$. Mais la seule chose qui compte dans sa
construction c'est la donn{\'e}e du demi-plan $P(e_1, e_2)$  bord{\'e} par la
droite 
$D_{\infty}(e_1, e_2)$ qui le contient. En fait, on peut partir de
n'importe quel 
demi-plan ferm{\'e} $P$ dont le bord $D_{\infty}^P$ est une droite
rationnelle. {\`A} l'int{\'e}rieur de ce demi-plan on consid{\`e}re la
droite $D^P$ parall{\`e}le {\`a} $D_{\infty}^P$  la plus proche de
$D_{\infty}^P$ qui contient des points de $N$. On consid{\`e}re alors sur
$D^P$ 
tous les couples de points successifs de $N$ : ce
sont des bases de $N$, qui permettent de construire le grand lotus
$\overline{\mathcal{L}}(P)$ associ{\'e} {\`a} $P$, comme union des lotus qui
leur correspondent. 
\medskip

Pour finir cette section, nous allons voir une relation entre le grand
lotus et les fractions continues de Hirzebruch-Jung. 

\begin{definition}
Consid{\'e}rons une suite $(v_0,v_1,...,v_{n+1})$ de
vecteurs de $N$, avec $n\geq 1$ . Cette suite est dite {\bf admissible} si les
conditions suivantes sont v{\'e}rifi{\'e}es :  

\begin{enumerate}
\item deux vecteurs successifs quelconques forment toujours une base de
$N$ ; 

\item  toutes ces bases d{\'e}finissent la m{\^e}me orientation de
$N_{\R}$ ; 

\item tous ces vecteurs sont contenus dans un m{\^e}me demi-plan ferm{\'e} 
bord{\'e} par la droite $\R v_0$. 
\end{enumerate}

\end{definition}

Les propri{\'e}t{\'e}s g{\'e}om{\'e}triques pr{\'e}c{\'e}dentes peuvent se
traduire num{\'e}riquement. Plus pr{\'e}cis{\'e}ment, on a la proposition
{\'e}l{\'e}mentaire suivante :

\begin{proposition} \label{vectrad}
 Pour chaque $ i \in \{1,2,3\}$, l'union des propri{\'e}t{\'e}s $1,  ..., i$
 est  {\'e}quivalente {\`a} l'union des propri{\'e}t{\'e}s $1', ..., i'$, o{\`u} :

 1'. $(v_0, v_1)$ est une base de $N$ et pour chaque $k \in
 \{1,...,n\}$, il existe $ \epsilon_k \in \{ 
 +1, -1\}$ tel que $v_{k-1} + \epsilon_k v_{k+1} = a_k v_k$, avec $a_k
 \in \Z$ ; 

 2'. de plus, $\epsilon_k=+1$ pour tous les $k \in \{1,...,n\}$ ; 

 3'. $a_k >0$ pour tout $k \in \{1,...,n\}$, d{\`e}s que $n\geq 2$ ;   
 $a_1\geq 0$ pour $n=1$ ;  $[a_1,..., a_k]^- >0$ pour tout $k \in
 \{1,...,n-1\}$ et $[a_1,...,a_n]^-\geq 0$. 

De plus, deux suites admissibles de vecteurs sont isomorphes par une
transformation lin{\'e}aire des r{\'e}seaux ambiants si et seulement si
les suites d'entiers associ{\'e}es co{\"\i}ncident. 
\end{proposition}

Ceci permet de parler de \emph{suites admissibles d'entiers}
$(a_1,...,a_n) \in \N^n$. La proposition pr{\'e}c{\'e}dente montre qu'une
suite admissible d'entiers est un invariant complet des suites
admissibles de vecteurs d'un r{\'e}seau bidimensionnel, {\`a}
isomorphismes de r{\'e}seaux pr{\`e}s.

Nous dirons  que la suite admissible $(v_0,...,v_{n+1})$ 
\emph{repr{\'e}sente $0$} si $v_1 + v_{n+1}=0$. Ceci est {\'e}quivalent au
fait que la suite repr{\'e}sentative $(a_1,...,a_n)$ v{\'e}rifie les
conditions 1', 2', 3' pr{\'e}c{\'e}dentes et que de plus :
  $$[a_1,...,a_n]^-=0 .$$

La proposition suivante fait le lien entre les suites admissibles
repr{\'e}sentant $0$ et les triangulations des polygones par des
diagonales. Elle peut se prouver par r{\'e}currence sur le nombre de
vecteurs, {\`a} l'aide du fait que les suites admissibles repr{\'e}sentant
$0$ s'obtiennent
{\`a} partir de la suite $(1,1)$ par un processus d'{\'e}clatements (voir
\cite[Appendice]{OW 77}).

\begin{proposition} \label{polycan}
  Soit $(v_0,...,v_{n+1})$ une suite admissible repr{\'e}sentant
  $0$ et soit $(a_1,...,a_n)$ la suite d'entiers associ{\'e}e. Notons
  par $P$ l'unique demi-plan ferm{\'e} contenant tous les 
  vecteurs de la suite. Soit $\mathrm{Pol}(v_0,...,v_{n+1})$ l'union
  du triangle $0v_1v_n$ et des
  p{\'e}tales du grand lotus $\overline{\mathcal{L}}(P)$ dont les
  int{\'e}rieurs intersectent l'un des 
  segments $]0, v_k]$, pour $k \in \{1,...,n\}$. C'est un polygone de
  sommets $0,v_1,...,v_n$ (dans cet ordre), triangul{\'e} par des
  p{\'e}tales du grand lotus,  et le nombre de
  p{\'e}tales arrivant au sommet $v_k$ est {\'e}gal {\`a} $a_k$, pour tout
  $k \in \{1,...,n\}$.
\end{proposition}

\begin{figure}[h!] 
\vspace*{5mm}
\labellist \small\hair 2pt 
  \pinlabel{0} at 255 -10
  \pinlabel{$v_0$} at 328 -10
  \pinlabel{$v_1$} at 367 42
  \pinlabel{$v_2$} at 406 128
  \pinlabel{$v_3$} at 288 74
  \pinlabel{$v_4$} at 224 75
  \pinlabel{$v_5$} at 98 166
  \pinlabel{$v_6$} at 120 120
  \pinlabel{$v_7$} at 9 145
  \pinlabel{$v_8$} at 130 40
  \pinlabel{$v_9$} at 184 -10
\endlabellist 
\centering 
\includegraphics[scale=0.50]{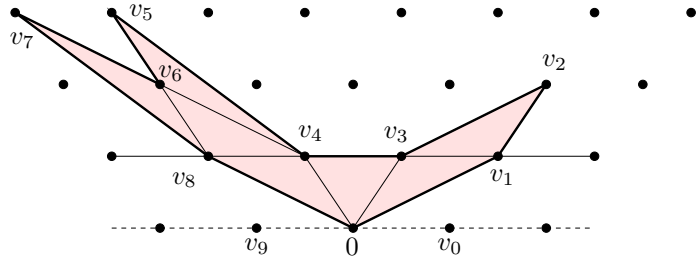} 
\vspace*{5mm} \caption{Polygone triangul{\'e} associ{\'e} {\`a} une suite
  admissible repr{\'e}sentant $0$} 
\label{fig:Triang} 
\end{figure}

\begin{example}
    Dans la Figure \ref{fig:Triang} est dessin{\'e} le polygone
    triangul{\'e} correspondant {\`a} la suite admissible repr{\'e}sentant
    $0$ suivante : $(2,1,3,4,1,3,1,3)$.
\end{example}

\begin{remark}
  La notion de suite admissible d'entiers a {\'e}t{\'e} introduite par Orlik
  et Wagreich \cite{OW 77}.  Celles qui de plus repr{\'e}sentent $0$
  sont r{\'e}apparues naturellement dans les travaux \cite{C 91} et \cite{S
    91} de Christophersen et Stevens sur la th{\'e}orie des
  d{\'e}formations des singularit{\'e}s quotient cycliques de
  surfaces. Ils montrent que ces suites correspondent bijectivement
  aux triangulations par des diagonales des polyg{\^o}nes orient{\'e}s
  ayant un point marqu{\'e} et ils utilisent la combinatoire de la
  triangulation associ{\'e}e dans leurs calculs. Par ailleurs, apr{\`e}s
  avoir {\'e}crit \cite{PP 07}, j'ai
  {\'e}t{\'e} amen{\'e} {\`a} me repr{\'e}senter ces suites par des suites de
  vecteurs, comme expliqu{\'e} dans la Proposition \ref{vectrad}. C'est
  {\`a} cette occasion que je me suis pos{\'e} la question de savoir s'il
  n'y avait pas moyen de construire canoniquement un polygone ayant la
  bonne triangulation {\`a} partir de la suite de vecteurs. 
\end{remark}

\section{Les lotus de dimension quelconque}
 \label{lotarb}

Dans cette section j'{\'e}tends la notion de lotus en dimension
quelconque, j'explique {\`a} quelle g{\'e}n{\'e}ralisation des fractions
continues cette construction donne lieu, et comment en d{\'e}duire une
mesure g{\'e}om{\'e}trique du contact des courbes monomiales.

\medskip
Partons d'une base (non-ordonn{\'e}e mais marqu{\'e}e par un ensemble $I$)
$\mathcal{B}:= (e_i)_{i \in I}$
d'un r{\'e}seau $N$ 
de rang $n\geq 2$, o{\`u} $I$ est un ensemble de cardinal
$n$. D{\'e}finissons les poly{\`e}dres convexes ferm{\'e}s suivants :

$\bullet$   $\sigma(\mathcal{B})$ : le c{\^o}ne strictement convexe de $N_{\R}$
engendr{\'e} par cette base ;

$\bullet$ $\Delta_0(\mathcal{B})$ : l'enveloppe convexe de l'ensemble
form{\'e} par $0$ et $\mathcal{B}$ ; c'est un simplexe $n$-dimensionnel ; 

$\bullet$ $\Delta(\mathcal{B})$ : l'enveloppe convexe de $\mathcal{B}$
; c'est un simplexe $(n-1)$-dimensionnel, unique facette de
$\Delta_0(\mathcal{B})$ ne contenant pas $0$ ;  

$\bullet$ $\Pi(\mathcal{B})$ : le parall{\'e}l{\'e}pip{\`e}de engendr{\'e} par
$\mathcal{B}$ ; ses sommets sont toutes les sommes de vecteurs
disjoints parmi les vecteurs de $\mathcal{B}$ ; 

$\bullet$  $\tau(\mathcal{B}) \: := \overline{\Pi(\mathcal{B}) \setminus
\Delta_0(\mathcal{B}) }$ ; c'est le \emph{p{\'e}tale} $n$-dimensionnel
associ{\'e} {\`a} la base $\mathcal{B}$ ; 

$\bullet$ $\phi(i_1,...,i_n)$ : pour chaque arrangement
$(i_1,...,i_n)$ des {\'e}l{\'e}ments de l'ensemble 
$I$, le simplexe $(n-1)$-dimensionnel dont
les sommets sont les {\'e}l{\'e}ments de la base : 
\begin{equation} \label{renum}  
\mathcal{B}(i_1,...,i_n) := (e_{i_1}, \: e_{i_1} + e_{i_2}, \:
  ... \: , \:  e_{i_1} + \cdots +   e_{i_n}).  
\end{equation}

Les simplexes $\phi(i_1,...,i_n)$  sont contenus dans le bord de
$\tau(\mathcal{B})$. Leur 
union avec $\Delta(\mathcal{B})$ constitue exactement la partie de
$\partial \tau(\mathcal{B})$ 
visible sans {\'e}crasement {\`a} partir de l'origine. C'est-{\`a}-dire que,
si on consid{\`e}re :
  $$ \psi: \sigma(\mathcal{B}) \setminus 0 \longrightarrow
  \Delta(\mathcal{B}), $$ 
la projection centrale sur le simplexe $\Delta(\mathcal{B})$ dont les 
sommets sont les points de $\mathcal{B}$, alors toutes les autres
faces maximales du bord du polytope $\tau(\mathcal{B})$ sont
{\'e}cras{\'e}es par $\psi$ en des polytopes de dimension strictement
inf{\'e}rieure. Les images par $\psi$ des simplexes $\phi(i_1,...,i_n)$
constituent exactement la subdivision barycentrique de
$\Delta(\mathcal{B})$.

\medskip 
\begin{figure}[h!] 
\vspace*{6mm}
\labellist \small\hair 2pt 
  \pinlabel{0} at 100 115
  \pinlabel{$e_1$} at 123 45
  \pinlabel{$e_2$} at 177 200
  \pinlabel{$e_3$} at 1 147
\endlabellist 
\centering 
\includegraphics[scale=0.60]{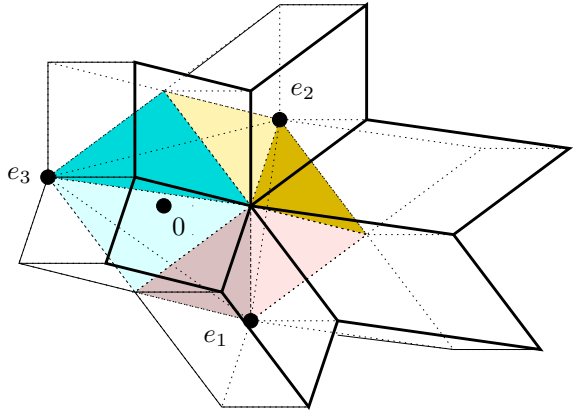} 
\vspace*{0mm} \caption{Les deux premi{\`e}res {\'e}tapes de
  construction du lotus tridimensionnel}    
\label{fig:Etages} 
\end{figure}

\medskip 
\begin{figure}[h!] 
\vspace*{6mm}
\labellist \small\hair 2pt 
  \pinlabel{$e_1$} at 450 -10
  \pinlabel{$e_2$} at -10 -10
  \pinlabel{$e_3$} at 220 450
  \pinlabel{$\Delta(e_1, e_2, e_3)$} at 380 250
\endlabellist 
\centering 
\includegraphics[scale=0.60]{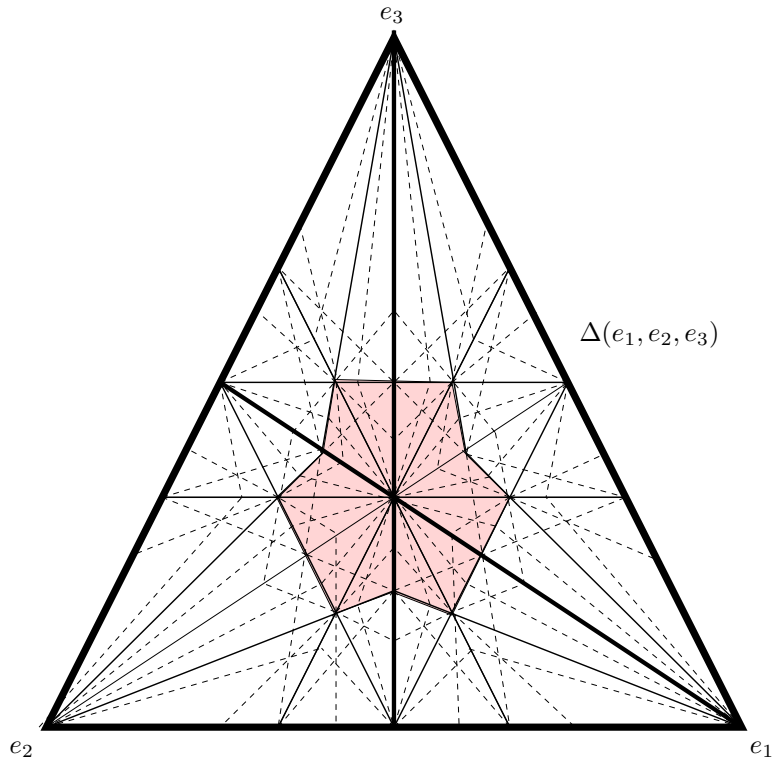} 
\vspace*{0mm} \caption{Le lotus tridimensionnel vu {\`a} partir de l'origine}   
\label{fig:Baryc} 
\end{figure} 

En partant de $\tau(\mathcal{B})$ et des nouvelles bases
$\mathcal{B}(i_1,...,i_n)$, on peut construire it{\'e}rativement
un complexe poly{\'e}dral infini contenu dans le c{\^o}ne
$\sigma(\mathcal{B})$. Plus pr{\'e}cis{\'e}ment, lors de la premi{\`e}re
{\'e}tape de la construction on construit $\tau(\mathcal{B})$. La
deuxi{\`e}me  {\'e}tape
de la construction consiste {\`a} rajouter tous les p{\'e}tales
$\tau(\mathcal{B}(i_1,...,i_n))$ {\`a} $\tau(\mathcal{B})$,
$(i_1,...,i_n)$ variant parmi les arrangements des {\'e}l{\'e}ments de $I$. Ces
p{\'e}tales `croissent' {\`a} partir des simplexes $\phi(i_1,...,i_n)$,
c'est pourquoi nous appelons ces derniers les \emph{simplexes de
  croissance} de 
la deuxi{\`e}me {\'e}tape de la construction. Remarquons que pour $n \geq
3$, ces simplexes de croissance ne sont pas des facettes du p{\'e}tale
$\tau(\mathcal{B})$, mais des demi-facettes. Ensuite on it{\`e}re...

\begin{definition}
Le complexe polyh{\'e}dral pr{\'e}c{\'e}dent est appel{\'e} {\bf le lotus} de
dimension $n$ 
associ{\'e} {\`a} la base $\mathcal{B}$ du r{\'e}seau $N$. On le notera
$\mathcal{L}(\mathcal{B})$. 
\end{definition}

\begin{example}
Dans la Figure \ref{fig:Etages}
sont dessin{\'e}s le p{\'e}tale initial et les 6 p{\'e}tales de la deuxi{\`e}me {\'e}tape de
construction. Sont colori{\'e}s les simplexes de croissance. 
Dans la Figure
\ref{fig:Baryc} sont dessin{\'e}s les p{\'e}tales des quatre premi{\`e}res {\'e}tapes de la
construction du lotus, tels que vus {\`a} partir de l'origine. C'est-{\`a}-dire que
sont dessin{\'e}es les images par $\psi$ des simplexes de croissance de
ces p{\'e}tales. Est colori{\'e}e la projection
de l'union des faces des p{\'e}tales de la deuxi{\`e}me {\'e}tape de
construction dont les 
ar{\^e}tes sont indiqu{\'e}es en traits gras dans la Figure
\ref{fig:Etages}. 
\end{example}

Le lotus $n$-dimensionnel permet de g{\'e}n{\'e}raliser en dimension
quelconque l'interpr{\'e}tation g{\'e}om{\'e}trique donn{\'e}e dans la
Proposition \ref{gainefrac}, 2) des fractions continues usuelles. En
effet, la notion de \emph{gaine} s'{\'e}tend en toutes dimensions :
  
\begin{definition}
   Un {\bf sous-p{\'e}tale} est une face lat{\'e}rale d'un p{\'e}tale,
   c'est-{\`a}-dire une face restreinte {\`a} laquelle $\psi$ n'est pas un
   hom{\'e}omorphisme sur son image. 
   Soit $l \subset N_{\R}$ une demi-droite d'origine $0$, contenue
   dans l'int{\'e}rieur du c{\^o}ne $\sigma(\mathcal{B})$. Sa {\bf gaine}
   $G(l)$ est par d{\'e}finition l'union des p{\'e}tales et sous-p{\'e}tales
   du  lotus $\mathcal{L}(\mathcal{B})$ dont les int{\'e}rieurs intersectent
   $l$. 
\end{definition}

\medskip 
\begin{figure}[h!] 
\vspace*{6mm}
\centering 
\labellist \small\hair 2pt 
  \pinlabel{$e_2$} at 0 -10
  \pinlabel{$e_1$} at 446 -10
  \pinlabel{$e_3$} at 220 455
  \pinlabel{$e_2 + e_3$} at 364 222
  \pinlabel{$e_1 + e_2 + e_3$} at -10 195
  \pinlabel{$e_1 + e_2 + 2 e_3$} at 70 366
  \pinlabel{$e_1 + 2 e_2 + 3 e_3$} at 354 355
  \pinlabel{$3 e_1 + 4 e_2 + 6 e_3$} at 443 173
\endlabellist 
\includegraphics[scale=0.50]{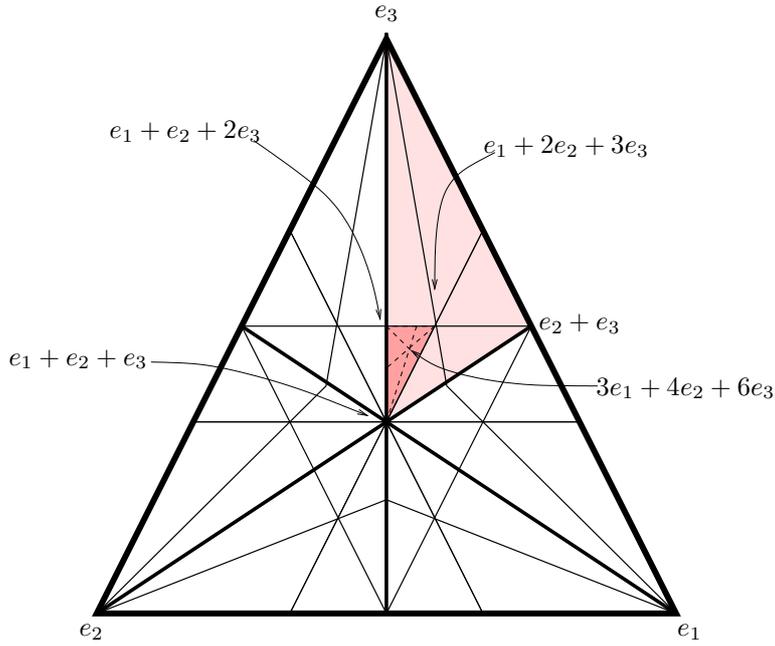} 
\vspace*{5mm} \caption{Une gaine tridimensionnelle vue {\`a} partir de
  l'origine}    
\label{fig:Gaine3} 
\end{figure}

On peut {\`a} nouveau associer {\`a} $G(l)$ la suite des sous-p{\'e}tales que l'on
rencontre en s'{\'e}loignant de $0$ le long de $l$. La suite des
dimensions de ces sous-p{\'e}tales est d{\'e}croissante. Pour $n=2$ elle
est constante, mais d{\`e}s $n=3$ elle ne l'est plus n{\'e}cessairement :
il y a une sous-suite initiale de p{\'e}tales de dimension $n$, suivie
d'une sous-suite de sous-p{\'e}tales de dimension $n-1$, puis une autre
de dimension $n-2$, etc. 

Notons par $\tau_0, \tau_1, \tau_2,...$ la suite des p{\'e}tales de
dimension $n$ de la gaine $G(l)$ et par $\Delta_0, \Delta_1,
\Delta_2,...$ la suite des simplexes de croissance associ{\'e}s. Gr{\^a}ce
{\`a} la formule (\ref{renum}), on peut param{\'e}trer cette suite par une
suite d'\emph{arrangements} $\alpha_0, \alpha_1, \alpha_2,...$ des
{\'e}l{\'e}ments de l'ensemble $I$, {\'e}tendant le fait que l'on
param{\`e}tre $\tau(\mathcal{B}(i_1,...,i_n))$ par l'arrangement
$(i_1,..., i_n)$. G{\'e}om{\'e}triquement, cela correspond au fait
qu'{\'e}tant donn{\'e} un simplexe de dimension $n$ de la subdivision
barycentrique de $\Delta(\mathcal{B})$, il admet un unique sommet en
commun avec $\Delta(\mathcal{B})$, et plus g{\'e}n{\'e}ralement, pour
chaque $k \in \{0,..., n-1\}$, une unique face de dimension $k$
contenue dans une face de dimension $k$ de
$\Delta(\mathcal{B})$. Ainsi, $\alpha_k$ est l'arrangement
correspondant au simplexe de croissance de $\tau_{k+1}$. Cette
param{\'e}trisation {\'e}tend en dimension plus grande l'association d'une
suite de symboles $\delta$ et $\gamma$ {\`a} une demi-droite $l \subset
\sigma(\mathcal{B})$ faite en dimension $2$ (voir le paragraphe qui
pr{\'e}c{\`e}de l'Exemple \ref{exgaine}). En effet, dans ce cas, si
$I=\{1,2\}$, le symbole $\delta$ correspond {\`a} l'arrangement $(1,2)$
et $\gamma$ {\`a} l'arrangement $(2,1)$.

\begin{example}
  Dans la Figure \ref{fig:Gaine3} est repr{\'e}sent{\'e}e la vue {\`a}
  partir de $0$ (c'est-{\`a}-dire son image par la projection $\psi$) de
  la gaine de la demi-droite rationnelle $\R_+(3e_1 + 4 e_2 + 6
  e_3)$. Cette gaine est compos{\'e}e de trois p{\'e}tales, param{\'e}tr{\'e}s par
  la suite $(3,2,1), (1,3,2), (1,3,2)$ d'arrangements des {\'e}l{\'e}ments
  de l'ensemble $\{1,2,3\}$. 
\end{example}

Chaque fois que la dimension des sous-p{\'e}tales constituant la gaine
baisse, on obtient des suites d'arrangements d'{\'e}l{\'e}ments d'un
ensemble de cardinal plus petit que celui qui pr{\'e}c{\'e}dait. De cette
mani{\`e}re, la notion de \emph{fraction continue} associ{\'e}e {\`a} une
demi-droite rationnelle de $\sigma(e_1,e_2)$ se retrouve remplac{\'e}e
en dimension plus grande par une notion \emph{d'arrangements
  continus}. 

\medskip
Pour finir, je voudrais mentionner une application de la notion de gaine {\`a} la
th{\'e}orie des singularit{\'e}s. Une \emph{courbe monomiale} de $\C^n$
est une courbe d{\'e}finie par une param{\'e}trisation de la forme $t
\rightarrow (t^{m_1},..., t^{m_n})$. Appelons $(m_1,...,m_n) \in \N^n$
son \emph{exposant}. On peut repr{\'e}senter g{\'e}om{\'e}triquement le
\emph{contact} de deux courbes 
monomiales {\`a} l'aide de la partie initiale commune des gaines de
leurs exposants. Ce contact peut aussi {\^e}tre repr{\'e}sent{\'e}
symboliquement en associant {\`a} la suite des p{\'e}tales de cette partie
initiale commune la suite correspondante d'arrangements.


\begin{thebibliography}{00} 
\bibitem{A 98} Arnold, V.I. \textit{Higher dimensional continued
    fractions.} Regular and chaotic dynamics \textbf{3}, 3 (1998),
    10-17. 

\bibitem{BK 86} Brieskorn, E.,  Kn{\"o}rrer, H. \textit{Plane algebraic
    curves.} Translated from the German by John Stillwell. Birkh{\"a}user
  Verlag, 1986.  

\bibitem{CC 05} Campillo, A., Castellanos, J. \textit{Curve
    singularities. An algebraic and geometric approach.} Hermann,
  2005.

\bibitem{CGL 92} Campillo, A., Gonzalez-Sprinberg, G., Lejeune-Jalabert,
M. \textit{Amas, id{\'e}aux {\`a} support fini et cha{\^\i}nes toriques.}
C. R. Acad. Sci. Paris S{\'e}r. I Math.  \textbf{315}  (1992),  no. 9, 987-990.  

\bibitem{CGL 96} Campillo, A., Gonzalez-Sprinberg, G., Lejeune-Jalabert,
M. \textit{Clusters of infinitely near points.}  Math. Ann. \textbf{306}
(1996),  no. 1, 169-194.  

\bibitem{CGM 09} Campillo, A., Gonzalez-Sprinberg, G., Monserrat,
  F. \textit{Configurations of infinitely near points.} A
  para{\^\i}tre dans S{\~a}o Paulo Journ. Math. Sciences.

\bibitem{CA 00} Casas-Alvero, E. \textit{Singularities of plane
    curves.} Cambridge Univ. Press, 2000.

\bibitem{C 91} Christophersen, J. A. \textit{On the components and
    discriminant of the versal base space of cyclic quotient
    singularities.} Dans  \textit{ Singularity theory and its
    applications}, Part I  (Coventry, 1988/1989), LNM \textbf{1462},
  Springer-Verlag, 1991, 81-92.

\bibitem{DV 36} Du Val, P. \textit{Reducible exceptional curves.}
  Amer. J. Math. \textbf{58} (1936), 285-289. 

\bibitem{EC 18} Enriques, F., Chisini, O. \textit{Lezioni sulla teoria
    geometrica delle equazioni e delle funzioni algebriche.} Libro II,
  Zanichelli, 1918. 

\bibitem{FJ 04} Favre, C., Jonsson, M. \textit{The valuative tree.}
  LNM \textbf{1853}. Springer-Verlag, Berlin, 2004. 

\bibitem{GP 09} Garc{\'\i}a Barroso, E., Popescu-Pampu, P. \textit{The
    kite of a plane curve singularity.} En pr{\'e}paration.

\bibitem{GKP 94} Graham, R. L., Knuth, D. E., Patashnik,
  O. \textit{Concrete Mathematics.} Second Edition, Addison-Wesley,
  1994.   

\bibitem{HT 85} Hatcher, A., Thurston, W. \textit{Incompressible
    surfaces in 2-bridge knot complements.} Invent. Math. \textbf{79}
  (1985), 225-246. 

\bibitem{H 73} Hironaka, H. \textit{La vo{\^u}te {\'e}toil{\'e}e.} Dans
  \textit{Singularit{\'e}s {\`a} Carg{\`e}se.} Ast{\'e}risque \textbf{7-8},
  SMF, 1973. 

\bibitem{H 63} Hirzebruch, F. \textit{The topology of normal
    singularities of an algebraic surface.} S{\'e}m. Bourbaki
  \textbf{250} (1962/63). 

\bibitem{H 00} Honda, K. \textit{On the classification of tight
    contact structures I.} Geometry \& Topology \textbf{4} (2000),
  309-368. 

\bibitem{K 96} Klein, F. \textit{{\"U}ber eine geometrische Auffassung
    der gew{\"o}hlischen Kettenbuchentwicklung.}
    Nachr. Ges. Wiss. G{\"o}ttingen. Math.-Phys. Kl. \textbf{3} (1895),
    357-359. French translation: \textit{Sur une repr{\'e}sentation
    g{\'e}om{\'e}trique du d{\'e}veloppement en fraction continue
    ordinaire.} Nouvelles Annales de Math{\'e}matiques (3) \textbf{15}
    (1896), 327-331.

\bibitem{LJ 95} Lejeune-Jalabert, M. \textit{Linear systems with
    infinitely near base conditions and complete ideals in dimension
    two.} Dans \textit{Singularity theory.} Trieste, 1991. D. T. L{\^e},
  K. Saito, B. Teissier eds., World Scientific, 1995, 345-369. 

\bibitem{M 61} Mumford, D. \textit{The topology of normal
    singularities of an algebraic surface and a criterion for
    simplicity.}  Inst. Hautes {\'E}tudes Sci. Publ. Math.  \textbf{9} (1961)
  5-22. 

\bibitem{N 75} Noether, M. \textit{{\"U}ber die 
singularen Werthsysteme einer algebraischen Function und die singularen 
Punkte einer algebraischen Curve.}  Math. Annalen \textbf{9} (1875), 166-182.


\bibitem{OW 77} Orlik, P., Wagreich, P. \textit{Algebraic surfaces
    with $k^*$ action.} Acta Math. \textbf{138} (1977), 43-81. 

\bibitem{PP 07} Popescu-Pampu, P. \textit{The geometry of continued
    fractions and the topology of 
     surface singularities}. Dans \textit{Singularities in Geometry
     and Topology 2004}. Advanced Studies in Pure Mathematics
   \textbf{46}, 2007, 119-195.

\bibitem{S 82} Series, C. \textit{Non-euclidean geometry, continued
    fractions and ergodic theory.} Math. Intelligencer \textbf{4}
  (1982) 1, 24-31. 

\bibitem{S 85} Series, C. \textit{The modular surface and continued
    fractions.} J. London Math. Soc. (2) \textbf{31} (1985), 69-80. 

\bibitem{S 91} Stevens, J. \textit{On the versal deformation of cyclic
    quotient singularities.} In \textit{Singularity theory and its
    applications.} Part 1 (Coventry, 1988/1989), LNM \textbf{1462},
  Springer-Verlag, 1991, 302-319. 

\bibitem{W 04} Wall, C. T. C. \textit{Singular points of plane
    curves.} London Mathematical Society Student Texts, \textbf{63}. 
   Cambridge University Press, 2004. 

\bibitem{Z 38} Zariski, O. \textit{Polynomial ideals defined by
    infinitely near base points.} Amer. J. Math. \textbf{60} (1938),
  151-204. 

\end{thebibliography}
\end{document}